\input amstex




\font\rm=cmr10 \rm

\font\bf=cmb10
\font\Rm=cmr9 at 11pt
\rm
\font\it=cmsl9 at 10pt
 at 7pt

\font\Rrm=cmr17 at 16pt
   \font\Rm=cmr12 at 11.5pt

\long\def\Pf{\par\noindent {\it Proof.} }
\def\({\left(}
\def\){\right)}
\def\st{such that }
\def\qed{\hfill$\bullet$\vskip 4pt}
\def\quotes#1{{\lq\lq #1\rq\rq}}
\def\brcs#1{\left\{ #1\right\}}

\def\eqv{\sim}

\def\wrt{with respect to }
\def\:{\,:}

\def\bb{\text{\bf b\,}}

\def\rr{{\Cal V}}
\def\EE{{\Cal E}}

\def\C{\text{\bf C}}

\def\Re{{\text{Re}\,}}
\def\Im{{\text{Im}\,}}

\def\R{\text{\bf R}}

\def\Arrow #1;#2.{#1\:#2 \to }

\def\Set#1#2{\brcs{#1 \left|\vphantom{#1 #2} \right.#2}}



\def\Rrr#1,#2{{\Cal J}_{#1,#2}}
\def\slfrac#1#2{{\raise -.07 ex\hbox{$^{#1}$}}\!/\raise .35 ex \hbox{${}_{#2}$}}
\def\ssf #1/#2{\slfrac {#1}{#2}}

\def\pd #1,#2.{\frac {\partial #1}{\partial #2}}

   \long\def\Lem
#1.#2\par{\vskip4pt{\baselineskip=13pt\font\it=cmsl12 at
11.5pt\Rm
   \noindent {\rm \uppercase{#1}} #2\vskip3pt

   }}

\long\def\Title #1\par {\noindent{\Rrm #1}\vskip 9pt}

 \long\def\SubT #1.{\noindent {\it #1\/} } 
 
 \long\def\SecT
#1\par{\vskip 3pt \noindent {\bf #1}\vglue1pt
   \noindent}

\long\def\subtitle #1.{\vskip 2pt \noindent {\it #1}}

\long\def\Rmk#1\par{\vskip 1pt \noindent {\it
Remark.} #1\vskip2pt}


\scrollmode\NoBlackBoxes
\magnification=1100
\long\def\Abstract #1\par%
{\vskip .2 true cm{\leftskip 1 true in \rightskip 1 true in \font\rm=cmr8 \rm
\baselineskip=1pt \font\it=cmsl8 \font\bf=cmb8
\parindent=0em {\bf Abstract} #1

}}
\comment
\font\rm=Times at 10pt

\font\bf=TimesB
\font\Rm=Times at 11pt
\rm
\font\it=TimesI at 10pt
\endcomment

\long\def\Pf{\par\noindent {\it Proof.} }
\def\({\left(}
\def\){\right)}
\def\st{such that }
\def\qed{\hfill$\bullet$\vskip 4pt}
\def\quotes#1{{\lq\lq #1\rq\rq}}
\def\brcs#1{\left\{ #1\right\}}
\def\Set#1#2{\brcs{#1 \left|\vphantom{#1 #2} \right.#2}}

\def\eqv{\sim}
\def\C{\text{\bf C}}

\def\Re{\text{Re\,}}

\def\Im{\text{Im\,}}

\def\wrt{with respect to }
\def\:{\,:}
\def\Arrow #1;#2.{#1\:#2 \to }


\def\R{\text{\bf R}}

\def\Rrr#1,#2{{\Cal J}_{#1,#2}}

\def\slfrac#1#2{{\raise -.07 ex\hbox{$^{#1}$}}\!/\raise .35 ex \hbox{${}_{#2}$}}
\def\ssf #1/#2{\slfrac {#1}{#2}}

\def\EE{{\Cal E}}
\def\pd #1,#2.{\frac {\partial #1}{\partial #2}}
\def\rr{{\Cal V}}


   \long\def\Title #1\par {\noindent{\Rrm #1}\vskip 9pt}
 \long\def\SubT #1.{\noindent {\it #1\/} }   \long\def\SecT
#1\par{\vskip 3pt \noindent {\bf #1}\vglue1pt
   \noindent}
\long\def\subtitle #1.{\vskip 2pt \noindent {\it #1}}

\long\def\Rmk#1\par{\vskip 1pt \noindent {\it
Remark.} #1\vskip2pt}


\def\One{1}
\def\oneone{\One.1}
\def\onetwo{\One.2}
\def\onethr{1.3}
\def\onefou{1.4}
\def\onefiv{1.5}
\def\onesix{1.6}
\def\onesev{1.7}

\def\throne{2.1}


\def\RRR{{\Cal R}}
\def\RRRRR{{\Cal S}}
\def\bb{{\text{{\bf b}}}}

\def\bb{\text{\bf b\,}}

\def\sign #1{\text{sign}\(#1\)}


\def\bb{\text{\bf b\,}}

\Title Arguments of zeros of highly log concave polynomials

\noindent {\it David Handelman%
\plainfootnote{\rm*}{\rm Supported by an NSERC Discovery Grant}}\vskip 4pt

{\leftskip = 1 true in
\rightskip =1 true in \parindent = 0em
{\it Abstract.} For a real polynomial $p = \sum_{i=0}^{n} c_ix^i$ with no negative coefficients and $n\geq 6$, let
$\beta (p) =
\inf_{i=1}^{n-1} c_i^2/c_{i+1}c_{i-1}$ (so $\beta (p) \geq 1$ entails that $p$ is log concave). If $\beta(p) >
1.45\dots$, then all roots of $p$ are in the left half plane, and moreover, there is a function $\beta_0 (\theta)$
(for $\pi/2 \leq \theta \leq \pi$) \st $\beta \geq \beta_0(\theta)$  entails all roots of $p$ have arguments in
the sector $| \arg z| \geq  \theta$ with the smallest possible $\theta$; we determine exactly what this function
(and its inverse) is (it turns out to be piecewise smooth, and quite tractible). This is a one-parameter extension
of Kurtz's theorem (which asserts that
$\beta
\geq 4$ entails all roots are real). We also prove a version of Kurtz's theorem with real (not necessarily
nonnegative) coefficients. {\par}MSC2010: 26C10, 30C15, 05E99

}\vskip 10pt 

\noindent As an outgrowth of a question concerning a class of analytic functions, we
 give criteria for all roots of real polynomials to lie in a sector of the
form $\Set{z \in \C}{|\arg z| > \theta}$, at least for $\pi \geq \theta
\geq \pi/2$ and asymptotically as $\theta \to 0$. The criteria depend only
on log concavity of the coefficients.

Specifically, if $f = \sum_{i= 0}^N c_i x^i$ (of degree $N \geq 6$) is a polynomial
with positive coefficients,
let $\beta := \inf_{i=1}^{N-1} c_i^2/c_{i+1} c_{i-1} $, and assume $\beta
> 1$. Then there is $\theta > 0$ \st for all roots, $z$, of $f$,   $|\arg
z| > \theta$ (where $\arg$ is the principal value, i.e., $\arg $ takes on
values in $(-\pi,\pi]$). The function $\beta \mapsto \theta$ is determined
exactly for $\pi/2 \leq \theta \leq \pi$. For example, if $\beta = 1 +
\sqrt 2$, then all roots of $f$ lie in the sector $|\arg z| > 3\pi/4$,
while if $\beta = 2$, then all roots lie in $|\arg z| > 2\pi/3$, and
moreover, these numbers are sharp.

We also show that if the  $c_i$ are assumed merely to be complex, then if
$\beta:= \inf |c_i|^2/|c_{i+1} c_{i-1}| \geq 4.45\dots$ (a root of a
transcendental equation), then $f$ has only simple roots and can be
located within specific annuli), and moreover, if the $c_i$ are real, then
all roots of $f$ are real. This is an extension of Kurtz's theorem, which
states that if the $c_j $ are all positive and $\beta > 4$, then all roots
are real. We also provide minor improvements on this result.

Then we consider in section 2 an old question [P] and [CC, section 4] (I am indebted to 
Tom Craven for these references). Form the entire
function (or the polynomial) $ g_{\beta} = \sum c_i x^i$ wherein the
quotients $\beta: = c_i^2/c_{i+1}c_{i-1}$ do not change in $i$. For what
values of $\beta$ does $g_{\beta}$ have only real roots? We provide an
answer, $\beta \geq \beta_0$, with $\beta_0  $ determined to 24
places (and show how to improve this), but unfortunately $\beta_0$ does
not appear to be connected to anything else. However, it does yield
apparent paradoxes; for example, there exists a polynomial (of any degree
exceeding $5$) $g$ for which $\beta (g) > 3.99$, but which has nonreal
roots; however, since $\beta_0 < 3.3$, any $g_{\beta}$ with $\beta > 3.3$
will have only real roots.

The original question that led to this article, was the determination of
conditions on a polynomial
$g$ guaranteeing all roots lie in the sector $|\arg z - \pi| < \pi/4$, and by the
result cited above, $\beta \geq 1 + \sqrt 2$ (in the presence of $N \geq
5$) is sufficient. This question itself emanated from a result in [H], guaranteeing
that a polynomial with $g(1) = 1$ belong to a class of analytic functions known there
as $\EE$, which play a role in classification criteria for AT ergodic transformations. 

\comment

While $\EE$ is closed under products (even infinite ones), it can happen
that the product of two non-members, $g$, $h$, satisfying $g(1) = h(1) =
1$ can belong to $\EE$. However, it is obvious that if $g$ (satisfying
$g(1) = 1$) does not belong, then neither does any nonzero power, positive
or negative). There is a  class of non-members which are relevant for the
classification results in [H], in fact, the definition of $\EE$ is used to
obtain approximate factorization results by this class.
For each integer $n$,  set $g_n = (1-x^n)/n(1-x)$ (the $n$ in the
denominator is the normalization, so that $g_n (1) = 1$). The degree of
$g_n$ is $n-1$ (an unfortunate consequence of the definitions). It is easy
to verify for each $n \geq 3$, that $g_n$ does not belong to $\EE$, and
therefore neither does any power.

Much more interesting, and the main result of section 3, is that no
product of any number of $g_n$s (with all $n \geq 3$, but no restriction
on multiplicities) belongs to $\EE$---despite the fact that all the $g_n$
are strongly unimodal and the sequences of coefficients of the products 
are very well behaved). In fact, we prove a stronger result, namely that
there exists $z_0$ on the unit circle \st for all $n \geq 3$,  $|g_n
(z_0)|^2 > \exp(-\rr(g_n)|1-z_0|^2)$, so the family of $g_n$s are
uniformly not in $\EE$, and therefore neither is any product of them. (We
establish this by showing the set of such $z_0$s has nonzero measure.)

\noindent Let $D$ be a neighbourhood of the unit circle, and let $\Arrow g;D.\C$ be $C^2$ at $1$ and satisfy $g(1)
= 1$. Define the {\it variance\/} of $g$ (more accurately, the variance of the
distribution whose moment generating function is $g$)
$$
\rr(g) = g''(1) + g'(1) - (g'(1)^2).
$$
Obviously if $g$ has only real coefficients in its Maclaurin expansion,
then $\rr(g)$ is real.

Define 
$$
\EE = \Set{\text{entire } \Arrow g;\C.\C}{g(1) = 1  \text{ and for all $z
\in S^1$, $|g(z)|^2 \leq \exp (-|1-z|^2 \Re \rr(g))$ } }. 
$$
This notion was introduced in [H]. All products of functions of the form $1-\alpha + \alpha x$ (where
$\alpha$ is any real number) and  $\exp (N(x-1))$ belong to $\EE$  (in
fact, $\EE$ is closed \wrt multiplication and \wrt uniform convergence on
compact subsets), but
$(1+ x + x^2)/3$ does not (nor does any of its powers). Moreover, if $q$
is a quadratic whose roots are on the unit circle and $q$ has strongly
unimodal coefficients, then with probability exceeding $.93$,  $q/q(1)$
will belong to $\EE$. Among real quadratics $q$ 
necessary  and sufficient conditions for $q/q(1)$ to belong to $\EE$ are
known [H; Appendix E].

We discuss two aspects of membership in $\EE$. The first section describes 
a criterion depending only locally on  the coefficients guaranteeing that 
all zeros of a Maclaurin series, $p$, lie in a sector of the form $|\arg z| > \theta$.
When $\theta \geq 3\pi/4$, in this context it implies further that $p/p(1)$ belongs
to $\EE$. In fact, we obtain much broader results about locations of zeros.

Suppose that $p = \sum_{j=0}^N c_j z^j$ ($6 \leq N \leq \infty$) is entire. 
If all $c_j$ are real and $\inf_{1 \leq j < N} |c_j^2/c_{j+1}c_{j+1}| \geq 4.45\dots$
(the actual number is the root of a transcendental equation), then all roots of $p$ 
are real. If we restrict to positive $c_j$, the corresponding result is known as Kurtz's
theorem [K], and the infimum of the ratios need only be $4$ (in fact, we offer minor
improvements on Kurtz's theorem, which was originally stated with $> 4$). In both cases,
it is easy to show that $p/p(1) $ belongs to $\EE$ (in the first case, we must assume
$p(1) \neq 0$).

We give sufficient conditions for all roots to be in the sector
$|\arg z | > \theta$ when $\pi/2 \leq \theta  <\pi$. For example, if $\theta = 3\pi/4$, and all $c_j > 0$,
then sufficient is that $N \geq 5$ and $\inf_{1 \leq j < N} c_j^2/c_{j+1}c_{j+1} \geq 1+\sqrt 2$
and this is sharp (that is, no smaller number than $1+\sqrt 2$ will do). It follows that if $p$ satisfies
this condition, then $p/p(1) $ belongs to $\EE$. We give an explicit formula for the optimal ratio, 
which is continuous, strictly monotone, and piecewise C${}^1$ as a function of $\theta$. We also give a roughly 
asymptotic formula as $\theta \to 0$. 

In section 2, we offer a surprising class of nonmembers of $\EE$. Let $g_n = (1-z^n)/n(1-z)$ (the $n$ in the
denominator is a normalization so that $g_n(1) =1$; it plays no other role). It is easy to check that $g_3 =
(1+z+z^2)/3$ is not in
$\EE$. We show that not only are none of
$\brcs{g_n}_{n\geq3}$ in $\EE$, but none of their products (allowing multiplicities) are as well. (Products of
normalized entire nonmembers of $\EE$ can belong to $\EE$.) We conclude with a couple of conjectures about
membership in  $\EE$.

In the Appendix, we discuss polynomials and entire functions for which $\beta = c_j^2/c_{j+1} c_{j-1}$ does not
change in $j$, and determine to within $10^{-24}$, the critical value of $\beta$, $B_0^2$, above which the
entire function and almost all the polynomials have only real zeros, and below which they have some nonreal
zeros.

\endcomment

\vskip 4pt
\noindent {\bf Section 1} Arguments of zeros
\vskip 2pt

\noindent 
\comment
Our initial aim in this section was to provide more examples of members
of
$\EE$ via factorization results, but one thing led to another \dots. A consequence of [H; Appendix E] is that if
$p$ is a polynomial all of whose zeros lie in the sector $\Set{z \in
\C}{|\arg z -
\pi| \leq \pi/4}$, then $p/p(1) $ belongs to $\EE$ (if we replace polynomial by entire
function of order at most one, the same result holds). 
\endcomment
Here we give
sufficient---but far from necessary---conditions for polynomials and
entire functions to have all their roots  in this sector, which however,
are easy to verify. However, we have more precise results for  sectors of the form
$\Set{z \in \C}{|\arg z| > \theta}$ for all $\theta$ with $\pi> \theta \geq \pi/2$.

A well-known theorem  due to Kurtz [K] asserts that if $p =
\sum_{i=0}^N c_i z^i$ is a polynomial of degree $N$ with $c_0 > 0$  and only nonnegative coefficients \st
for all $i=1, \dots, N-1$, the numbers $\beta_i(p):=
c_i^2/c_{i+1}c_{i-1}$ all exceed $4$, then all roots of $p$ are real (and
thus negative). This is extended to give similar type conditions (on the
ratios, $\beta_i(p)$) to guarantee that all the zeros lie in a sector of
the form $|\arg z - \pi| < \psi$ for $ 0 < \psi \leq \pi/2$. As a special
case, we show that if $N \geq 5$ and $c_i^2/c_{i+1}c_{i-1} \geq 1 + \sqrt
2$, then all zeros satisfy $|\arg z - \pi | < \pi/4$, so that the
corresponding $p/p(1)$ belongs to $\EE$. The number $1+ \sqrt 2$ is sharp
in the sense that for all $\epsilon > 0$, there exists a polynomial, $p$,
of degree $N$ with some roots  outside the  sector, yet with $c_i^2
/c_{i+1}c_{i-1} > 1 + \sqrt 2-\epsilon$ for all $i = 1,\dots, N-1$.

For $\theta = \pi - \psi$, define two functions,
$$\eqalign{
\RRR(\theta) & = \cases \text{largest positive real root of $X^2 - 2 \cos
\theta \cdot X^{3/2} + 2 \cos 2\theta = 0$}& \text{if one exists} \\ 1 &
\text{otherwise}
\endcases \cr
\RRRRR(\theta) & = \cases \text{largest positive real root of $X^3 +
\frac{\cos 3\theta/2}{\cos \theta/2}\cdot X^{2} + \frac{\cos
5\theta/2}{\cos
\theta/2} = 0$}& \text{if one exists} \\ 1 & \text{otherwise.}
\endcases \cr 
}$$ 
We will show that if $N \geq 6$ (or $N \geq 5$ if
$\theta$ is not too close to $\pi/2$), $\pi/2 \leq \theta < \pi$ and
$$
\inf_{i=1}^{N-1} \frac {c_i^2}{c_{i+1}c_{i-1}} \geq \max \brcs{4 \cos^2
\theta, 1 - 2 \cos \theta, \RRR(\theta), \RRRRR(\theta)},
$$ 
then all
roots of $p = \sum c_i z^i$ lie in the sector $|\arg z| < \theta$, and
moreover, this is sharp (in the sense of the example  with $\theta
= 3\pi/4$ and $1 -2\cos \theta = 1 + \sqrt 2$); this criterion can be slightly simplified, as in the statement of
Theorem \onesix.

For example, with $\theta = \pi/2$, all roots of $p$ lie in left half
plane if $\inf_i c_i^2/c_{i+1}c_{i-1}$ is at least as large as the
positive real root of $X^3 - X^2 -1 =0$, about $1.45\dots$. For $\theta
= 2\pi/3$, the corresponding lower bound is $2$.

These results extend in a very routine way to entire functions. When $g =
\sum c_i z^i$ with $c_i > 0$ is entire and $\inf c_i^2/c_{i+1}c_{i-1} >
1$, then $g$ is of order zero, and admits a factorization of the form
$g/g(0) = \prod (1 - z/z_k)$ where $z_k$ runs over the zeros.

Since we deal almost exclusively with polynomials all of whose zeros lie
in $|\arg z| \geq \pi/2$, these polynomials will automatically have no
negative coefficients, and if $|\arg z | \geq 2\pi/3$, the polynomials
will be strongly unimodal (that is, the sequence consisting of their coefficients will
be log concave). We give a result, Theorem
\oneone, along the lines indicated here that does not require the coefficients to be
nonnegative, but merely real, with the conclusion that the roots are all
real. Of course, there is a vast literature on polynomials and entire
functions all of whose zeros are real, but I couldn't find this
particular result in the literature (which of course does not mean that
it does not exist therein).

 The
lower bound given in Theorem \oneone, $\beta_0$, is likely not sharp, but
on the other hand, as was noted in [K], whatever the optimal value is, it
must be at least  $25/6 > 4$, because the polynomial $x^3 - 5x^2 + 6x +1$
has nonreal roots.

Define the function $F(r,\beta) = 1 + r + r^2/\beta + r^3/\beta^3 + \dots
=
\sum_{j=0}^{\infty} r^j/\beta^{j(j-1)/2}$; typically, $\beta > 2$ and $0
\leq r \leq 1$. It is perhaps accidental that the function $F_{\beta} (z)
= F(z,\beta)$ is an entire function (when $\beta > 1$) satisfying
$\beta_k = \beta^2$ (so if $\beta \geq 2$, then all of its roots are
real).

Suppose that $\brcs{D_k}$ is a summable sequence of positive real
numbers   and $\beta_j:= D_j^2/D_{j+1}D_{j-1} \geq \beta$ for some number
$\beta > 1$. Since the sequence is strongly unimodal, it is unimodal, and
let $m$ be the smallest mode, and let $k \geq m$, so that the sequence
$\brcs{D_j}_{j\geq k}$ is monotone nonincreasing. Set  $r = D_{k+1}/D_k$;
necessarily, $r \leq 1$. For $l > k+1$,
$$\eqalign{
\frac{D_l}{D_{l+1}} &= \frac{D_{l-1}}{D_l} \frac 1{\beta_{l-1}}, \qquad
\text{which iterates to} \cr & = \frac{D_{k+1}}{D_k} \frac
1{\beta_{l-1}\beta_{l-2}\dots \beta_{k+1}} = \frac
{r}{\prod_{j=k+1}^{l-1} \beta_j} .\cr }$$ If we also assume that $\beta_j
\geq \beta$ for all $j$, then we have $D_l/D_{l+1} \leq r\beta^{l-k-1}$.
Now $D_l/D_k$ telescopes as a product $(D_l/D_{l-1})(D_{l-1})/D_{l-2}
\cdots$, and we thus obtain that $D_l \leq r^{l-k-1} 
\beta^{(l-k-1)(l-k-2)}D_k$. Hence the mass of the tail, that is, the sum
$\sum_{j \geq k} D_j$ is bounded above by $F(r,\beta) D_k$.

Similarly, if $k \leq m$, the sequence is increasing, and on setting $s =
D_{k-1}/D_k$, we have $\sum_{j \leq k} D_j < F(s,\beta) D_k$; the
inequality is strict, because the sequence is finite to the left of $m$.

Now we are in position to prove a result via  an easy application of
RouchŽ's theorem, where we use a single monomial as the function to which
we compare the zeros. To avoid cluttering the statement of the theorem even
more than it currently is, we define
$$\eqalign{
\rho_k&:=\max\brcs{\frac{d_{k-2}}{rd_{k-1}}, \frac{d_k} {2 d_{k+1}F(r,\beta)} \(1- \sqrt{1 -
\frac{4F(r,\beta)}{\beta_k}}\)}\cr
R_k&:=  \min\brcs{\frac{rd_{k+1}}{d_{k+2}}, \frac{d_k} {2 d_{k+1}F(r,\beta)} \(1+ \sqrt{1 -
\frac{4F(r,\beta)}{\beta_k}}\)}, 
}$$
where $r = \beta^{-3/2}$.

\Lem Theorem \oneone. Let $f = \sum_{j = 0}^N c_j z^j$ (with complex
$c_j$) be an  entire function with  $N$ in $\brcs{3,4,\dots} \cup
\brcs{\infty}$. Suppose that  the sequence $(d_j =
|c_j|)$ satisfies
$\beta_j := d_j^2/d_{j+1}d_{j-1} \geq \beta_0$ for all $j$, where
$\beta_0 = 4.448505576\dots$ is the unique root of
$F(\beta^{-3/2},\beta)^2 = \beta$. Then all zeros of $f$ are simple, and
exactly one appears in the annulus $R_k < |z| < \rho_k$, and there are no
others. If additionally,
$c_j$ are all real, then all zeros of $f$ are real.

\Pf We will apply RouchŽ's theorem with $g(z) = c_k z^k$ (once for each
$k$); we will show that $|(f-g) (z)| < |g(z)|$ on the circle $|z| = R$,
with $R$ to be chosen appropriately. It follows that $f$ has exactly $k$
zeros in the disk, and by increasing $k$ to $k+1$, we will obtain a
larger disk containing exactly $k+1$ zeros, so there  must be exactly one
zero in the set-theoretic difference of the disks, that is, an annulus.

For unspecified $R > 0$, set $D_j = d_j R^j = |c_j z^j|$  (on $|z| = R$).
Fix $k$, and suppose that $d_k \geq d_{k+1} R, d_{k-1} /R$; that is, $k$
is a mode of the sequence $\brcs{d_j R^j}$. Define $r =
\max\brcs{D_{k+2}/D_{k+1}, D_{k-2}/D_{k-1}}$, so that $r \leq 1$. Set
$\beta = \inf \brcs{\beta_j}$. Then we have
$$\eqalign{
\sum_{j \geq 1} D_{k+j} & \leq D_{k+1} F(r, \beta)\cr
\sum_{j \geq 1} D_{k-j} & <  D_{k-1} F(r, \beta). \cr
}$$ 
Obviously, for
$z$ on the circle, $|g(z)| = d_k R^k = D_k$; equally obviously,
$|(f-g)(z)|
\leq \sum_{j \neq k} D_j$, so to show $|(f-g)(z)| <|g(z)|$, it would
suffice to show that $ D_k \geq F(r,\beta) (D_{k+1} + D_{k-1})$. We now
derive conditions on $R$, $r$, and $\beta$ to guarantee this; it is
equivalent to the quadratic inequality,
$$\eqalign{ 
R^2 d_{k+1} -&\frac{Rd_k}{F(r,\beta)} + d_{k-1} \leq 0
\quad{\text{that is,}}\cr
&\cases
\beta_k \geq 4F(r,\beta) &\text{ \hskip -4 true in and}\\
\frac{d_k}{2d_{k+1} F(r,\beta)}\(1 - \sqrt{1-\frac {4
F(r,\beta)}{\beta_k}}\) \leq R \leq
\frac{d_k}{2d_{k+1} F(r,\beta)}\(1 + \sqrt{1-\frac {4
F(r,\beta)}{\beta_k}}\). & \\
\endcases \cr 
}$$ 

To summarize, at this point, we require conditions on $R$ (which is to be
determined), and $r$ (also to be determined), as follows:
$$\eqalign{
\frac{d_{k-1}}{d_k} & \leq R \leq \frac{d_k }{d_{k+1}}\cr
\frac{d_k}{2d_{k+1}F(r\beta)} \(1 -
\sqrt{1-\frac{4F(r,{\beta_k})}{\beta_k}} \) & \leq R \leq
\frac{d_k}{2d_{k+1}F(r\beta)} \(1 +
\sqrt{1-\frac{4F(r,\beta_{k})}{\beta_k}} \) \cr
\frac{D_{k+2}}{D_{k+1}} , \frac{D_{k-2}}{D_{k-1} }&\leq r, \cr
}$$
and to make $r$ as small as possible subject to these conditions. The
third condition is equivalent to
$$
   \frac{d_{k-2}}{rd_{k-1}}  \leq R \leq  \frac{r d_{k+1}}{d_{k+2}}, \tag*
$$
which actually subsumes the first. Necessary and sufficient for the
existence of $R> 0$ satisfying (*) is simply $d_{k+1}d_{k-1}/d_{k+2}d_{k-2}
\geq 1/r^2$, that is, $ \beta_{k-1} \beta_k \beta_{k+1} \geq 1/r^2$. (The
left side will be at least $64$, so this allows $r$ at this stage to be
fairly small.) For the second condition and (*) to be compatible, we
require  the two inequalities (which now constrain $r$).
$$\eqalign{
  \frac{d_{k-2}}{rd_{k-1}} & \leq \frac{d_k}{2d_{k+1}F(r,\beta)} \(1 +
\sqrt{1-\frac{4F(r,\beta_{k})}{\beta_k}} \) \quad \text{and} \cr
\frac{d_k}{2d_{k+1}F(r,\beta)} \(1 -
\sqrt{1-\frac{4F(r,\beta_{\beta_k})}{\beta_k}} \) & \leq
 \frac{r d_{k+1}}{d_{k+2}}.\cr
}$$
These two inequalities can be written in the form
$$\eqalign{
\frac 1r & \leq \frac{\beta_k\beta_{k-1}}{2F(r,\beta_k)} \(1 + \sqrt{1 - \frac{4F(r,\beta_k)}{\beta_k}}\)\cr
\frac 1r & \leq \frac{2\beta_k F(r,\beta_k)} {1 - \sqrt{1 - \frac{4F(r,\beta_k)}{\beta_k}}} = 
\frac{2\beta_k F(r,\beta_k)}{4F(r,\beta_k)/\beta_k}\(1 + \sqrt{1 - \frac{4F(r,\beta_k)}{\beta_k}}\)\cr
& = \frac{\beta_k^2 F(r,\beta_k)}2\(1 + \sqrt{1 - \frac{4F(r,\beta_k)}{\beta_k}}\)\cr
}$$

Now suppose that $\beta \leq \beta_k, \beta_{k\pm1}$. Since $F(r, \beta) > 1$ (when $r > 0$), sufficient for both
these inequalities to hold is 
$$
\frac 1r \leq \frac{\beta^2}{2 F(r,\beta)}\(1 + \sqrt{1 - \frac{4F(r,\beta)}{\beta}}\).
$$
Therefore, it is sufficient to find $r$ and $\beta$ so that (for suitable $\beta$), the following hold:
$$\eqalign{
4 F(r,\beta) & \leq \beta\cr
\frac{1}r & \leq \min \brcs{\beta^{3/2},\frac{\beta^2 F(r,\beta)}2\(1 + \sqrt{1 -
\frac{4F(r,\beta)}{\beta}}\) } 
}$$

Normally, this would be hopeless; however, there is a trick, obtained by setting the two terms in the minimum to
each other. Let $\beta_0$ be the unique positive solution to $\beta = (F(\beta^{-3/2},\beta) + 1)^2$; a back of the
envelope calculation yields easily that $4.3 < \beta_0 < 4.5$. {\it Maple\/} yields $4.448505576\dots$ (convergence
is extremely fast). Setting
$r =
\beta_0^{-3/2}$, we see fairly quickly that all the relevant inequalities hold.

Hence if each of $\beta_k, \beta_{k\pm 1}$ are at least as large as $\beta_0$, then RouchŽ's theorem applies, and
we deduce that with 
$$\eqalign{
\max\brcs{\frac{d_{k-2}}{rd_{k-1}}, \frac{d_k} {2 d_{k+1}F(r,\beta)} \(1- \sqrt{1 -
\frac{4F(r,\beta)}{\beta_k}}\)} \leq R \leq\cr
&  \text{\hglue -1.5 true in} \leq \min\brcs{\frac{rd_{k+1}}{d_{k+2}}, \frac{d_k} {2 d_{k+1}F(r,\beta)} \(1+ \sqrt{1 -
\frac{4F(r,\beta)}{\beta_k}}\)}, 
}$$
the interval is nonempty (it could be a singleton), and $f$ has exactly $k$ zeros in $|z| < R$ and none on $|z| = R$.
Then $\rho_k$ is the left endpoint   and $R_k$ is the right endpoint.
Assume that
$\beta_j \geq \beta_0$ for all $j$. Necessarily $\rho_{k+1} > R_k$. It follows that on the annulus $\rho_k \leq |z|
\leq R_k$, $f$ has no zeros, and on the annulus $R_k < | z| < \rho_{k+1}$, $f$ has exactly one root. In particular,
all roots are simple. Moreover, if we additionally assume that all $c_j$ are real, then the roots must be real
(since nonreal roots come in conjugate pairs).
\qed

Since the argument is based on RouchŽ's theorem, the constant is unlikely
to be optimal.

Now we can slightly extend Kurtz's theorem [K]. By restricting to polynomials
of higher degree (at least $3$), we can replace the strict inequalities
that appear in the original statement by greater than or equal signs. We 
also give rather crude ranges for the locations of the zeros.

For a nonzero real number $r$, define $\sign r = -1$ if $r < 0$
and $\sign r = 1$ if $r > 0$.

\Lem Theorem \onetwo. Suppose that $f = \sum_{j =0}^N c_j z^j$ with $N \in
\brcs{3,4,\dots,} \cup \brcs{\infty}$ is entire, all $c_j > 0$, and for
all $j= 1, 2, \dots, N-1$, we have
$\beta_j:= c_j^2/c_{j+1}c_{j-1} \geq 4$. Then all roots of $f$ are simple
and real, and moreover, if $-x_k$ is the $k$th smallest (negative) root,
then
$$
\frac{c_{k-1}}{2c_k} \( 1 + \sqrt{1- \frac{4}{\beta_{k-1}}}\) \leq x_k
\leq
\frac{c_{k}}{2c_{k+1}} \( 1 - \sqrt{1- \frac{4}{\beta_{k}}}\),
$$ and this interval is nontrivial.

\Pf Define the closed interval
$$ J_k = \left[\frac{c_k}{c_{k+1}}\(1 - \sqrt{1 - \frac 4 {\beta_k}}\) ,
\frac{c_k}{c_{k+1}}\(1 + \sqrt{1 - \frac 4 {\beta_k}}\)\right].
$$ This is a singleton when $\beta_k =4$, otherwise it has nonzero
length. We will show that for $x$ in $J_k$, $\sign {f(-x)} = (-1)^k$.

By strong unimodality of $\brcs{c_j}_{j=0}^N$, the sequence $\brcs{x^j
c_j}$ is unimodal for any $x > 0$. Suppose positive  $x$ satisfies
$$
    c_k x^k \geq c_{k+1}x^{k+1} + c_{k-1} x^{k-1}; \tag *
$$ 
then the sequence $\brcs{x^j c_j}$ has a maximum at $j =k$, hence is
decreasing  for $j > k$ and increasing for $j < k$. Since the expansion
of $f(-x)$ is alternating, it follows immediately  that $|f(-x)
|$ is at least as large as  $\sum_{j\leq k-2} c_j(-x)^j + \sum_{j\geq
k+2} c_j(-x)^j$; at least one of these partial sums is not zero, since
degree $f$ is $N \geq 3$, and whenever one of the  partial sums or
$(-1)^k(c_k x^k - c_{k-1}x^{k-1} - c_{k+1}x^{k+1})$ is not zero, the sign is
simply $(-1)^k$. Hence $f(-x)$ would have sign $(-1)^k$.

Now (*) is equivalent to the quadratic inequality $x^2 - x c_k/c_{k+1} +
c_{k-1}/c_{k+1} \leq 0$; this in turn is equivalent to $(x-
c_k/2c_{k+1})^2 \leq (c_k^2 - 4 c_{k-1} c_{k+1})/4 c_{k+1}^2$. This has a
real solution if and only if $\beta_k \geq 4$, and in that case, the
solutions are precisely the points of $J_k$.

For $k = 0$, the corresponding result is the obvious $f(-x) > 0$ if $0 <
x \leq c_0/c_1$. Now it follows (easily) from strong unimodality of the
sequence, that the right endpoint of $J_k$  is less than the left
endpoint of $J_{k+1}$. Hence $f$ has at least $N+1$ sign changes on the
negative reals. Now suppose that $N$ is finite. It has at least $N$
distinct negative roots. Since the degree of $f$ is $N$, this must
account for all of them. Moreover, the roots can only occur between the
$J_k$, that is, between the right endpoint of $J_k$ and the left endpoint
of $J_{k+1}$, which yields the range in the statement of the theorem.

Now suppose that $N =\infty$ and $f$ is entire. Set $f_n = \sum_{j=0}^N
c_j z^j$; then $f_n \to f$ uniformly on compact subsets of $\C$, and of
course $f$ is not identically zero. Let $z_0$ be a root of $f$. If $z_0$
is not on the negative real axis, then it has  a neighbourhood which
misses then negative reals. Thus none of the $f_n$ have zeros on this
neighbourhood, hence $f$ cannot have a zero therein (since $f$ is not
identically zero), a contradiction. So all the real zeros of $f$ lie on
the negative reals, and we also know that the sign of $f(-x)$ does not
change on any $J_k$. So all zeros lie in the indicated sets, and there is
at least one (because of the sign changes) in each one. It remains to
show there can be no more than one.

This is again a consequence of uniform convergence on compact sets; let
$D_k$ be the open disk centred at the midpoint between $-J_k$ and
$-J_{k+1}$ and whose diameter joints the right endpoint of one to the
left endpoint of the other. On the bounding circle, note that $f_n'/f_n$
converges uniformly to $f'/f$, so the number of zeros enclosed also
converges; hence $f$ has just one zero in $D_k$.
\qed

Now we want to obtain conditions on the ratios ($c_i^2/c_{i+1}c_{i-1}$)
to guarantee that all the zeros are in sectors of the form $|\arg z | >
\theta$ for given $\theta$ in $[\pi/2,\pi)$. Let
$P_N$ denote the collection of real polynomials of degree
$N$ or less, topologized by identifying the polynomial $\sum c_j z^j$
with the point $(c_j)$ in $\R^{N+1}$.

\Lem Lemma \onethr. Suppose that $U$ is a nonempty open subset of $P_N$,
and
$W$ is a nonempty regular open (equal to the interior of its closure)
subset of $\C$ with the following properties.{

}\item{(a)} $U$ is connected and all elements of $U$ are degree $N$
\item{(b)} no element of $U$ has a zero on $\partial W$, the boundary of $W$
\item{(c)} there exists $p$ in $U$ \st all  of its zeros  lie in $W$.{

}\noindent Then all zeros of all members of $U$ lie in $W$.

\Pf Set $A = \Set{g \in U}{\text{all zeros of $g$ lie in $W$}}$; then $A$
is nonempty. We show both $A$ and $U\setminus A$ are open. For $g$ in
$A$, there exist tiny disks centred about its zeros, all of which,
including their closure, are contained in $W$. On each such disk, we may
assume that $g$ does not vanish on the bounding circle; say $\delta $ is
the infimum of the values of $|g|$ on the union of the circles. If $h$ is
in $U$ and $\left\|h - g\right\|  $ (the norm is the absolute sum of the
coefficients) is sufficiently small, then we can apply RouchŽ's theorem
to $h$ and $g$ (on the union of the disks), and so deduce that
$h$ has $N$ zeros within the union of the disks; since $h$ has degree $N$, this accounts for all of its zeros. Hence $h$
belongs to $A$.

The argument for $U\setminus A$ is similar, but we only have to work with one zero, $z_0$. If $g$ is in $U\setminus
A$ and $z_0$ is a zero of $g$ not in
$W$, then $z_0$ must lie in the complement of the closure of $W$ (since $W$ is regular, and $h$ has no zeros on $\partial
W$). Hence there exists a disk therein centred at $z_0$. The same RouchŽ's theorem argument (this time for a single
disk) yields that any $h$ sufficiently close to $g$ (in the coefficientwise norm) must have a zero in the disk. 

Since $U$ is connected and $A$ is open, $U = A$.
\qed

Let $\bb  := (b_1, \dots, b_{N-1})$ be a sequence of positive numbers, and define the following sets
$$\eqalign{
U_0 (\bb) & = 
\Set{c = (c_j)_{j=0}^N  \in (\R^{++})^{N+1}}
{\beta_j(c):=  \frac{c_j^2}{c_{j+1}c_{j-1}} > b_j, \ j = 1, 2,
\dots, N-1} \cr
C(\bb) & = \Set{(x_j)_{j=1}^{N-1}  \in (\R^{++})^{N-1}}{x_j > b_j, \ j = 1, 2,
\dots, N-1}.\cr
 }$$

Obviously $C(\bb)$ is  $(b_j, \infty)^{N-1}$ and $U_0(\bb)$ is open in
$\R^{N+1}$. Define $\Arrow \phi \equiv \phi(\bb); U_0 (\bb). (\R^{++})^2
\times C(\bb)$  via $\phi ((c_0, c_1, c_2, \dots, c_N)) = (c_0, c_1, \beta_1(c), \beta_2 (c), 
\dots, \beta_{N-1}(c))$.
This map is obviously well-defined and continuous. To construct its inverse, we note the recursive equations, $c_{j+1} =
c_j^2/\beta_j c_{j-1}$. Iterating these yields each $c_j$ ($j\geq 2$) as a function of $c_0$, $c_1$ and the $\beta_j$,
and all the denominators that appear are strictly positive. It is immediate that this yields a (trivially) rational
mapping (hence continuous) that is the inverse of $\phi$. In particular, $U_0 (\bb)$ is homeomorphic to $\R^{N+1}$, and
is thus connected!

Let $U(\bb)$ be the image of $U_0(\bb)$ in $P_N$, that is, associate to $c$ in $U_0$ the polynomial $f_c = \sum c_j
z^j$. So we have that $U(\bb)$ is an open and connected subset of $P(\bb)$.  Supppose $U(\bb)$ and $W$ satisfy the
conditions of the lemma. For each choice of   $c_0, c_1 > 0$, take the closure of the points in $U$ whose first two
coordinates  are $(c_0,c_1)$; then take the union over all strictly positive choices of $(c_0,c_1)$. The resulting set,
call it $V(\bb)$ is easily described: it is the set of points $c= (c_j)$ for which the corresponding $\beta_j \geq b_i$. 

Suppose that no members of $V(\bb)$ have zeros on $\partial W$; then it is easy to show that all zeros of all members of
$V(\bb)$ lie in $W$ (this is a bit surprising, since the latter is open). For if $z_0$ is a zero of $f$ in $V(\bb)$,
then $z_0$ cannot belong to $\partial W$ by hypothesis, so if $z_0 $ is
not in
$W$, it must lie in the complement of the closure of $W$, hence there is
a disk centred at it that lies entirely in the complement of the closure.
There exist $f_n$ in
$U(\bb)$ converging coordinatewise to $f$ (with the first two coordinates fixed), hence $f_n \to f$ on compact sets, and
once again, so $f$ can have no zeros in the disk, a contradiction. 

We state this as a corollary.

\Lem Corollary \onefou. Let $\bb = (b_1, \dots, b_{N_1})$ be an $N$-tuple
of real numbers with $b_i >1$, and set 
$$
V(\bb) = \Set{f = \sum_{j =  0}^N c_j z^j \in P_N}{c_j > 0, \
\frac{c_j^2}{c_{j+1}c_{j-1}} \geq b_i}.
$$ 
Suppose that $W$ is a regular
open subset of $\C$ with boundary $\partial W$, and there exists $f$ in
$V(\bb)$ that has all of its zeros in $W$. Suppose that every
$f$ in $V(\bb)$ has no zeros in $\partial W$. Then all the zeros of every $f$ in $V(\bb)$ lie in $W$.

For an angle $0 < \theta < \pi$, let
$W_{\theta}$ denote the open sector in $\C$, $\Set{z \in \C}{|\arg z| >
\theta}$ (the branch of $\arg z$ has values $-\pi  < \arg z \leq \pi$),
and let $L_{\theta}$ denote the ray in the upper half-plane,
$\Set{\lambda e^{i\theta}}{\lambda \geq 0}$. Then $\partial W_{\theta} = L_{\theta} \cup
\overline{L_{\theta}}$. Since the polynomials in $U(\bb)$ are real, to verify condition (b) in
Lemma \onethr, we need only verify that all polynomials therein have no
zeros on
$L_{\theta}$.

Now to verify $f$ has no zeros on $L_{\theta}$ for any $f$ in $U$, we can
make a further reduction. The reparameterization maps, $f \mapsto
f_{\lambda}$ (for each $\lambda > 0$), where $f_{\lambda}(z) = f(\lambda
z)$, do not change the
$\beta_j$ values. Hence $f$ belongs to $U(\bb)$ if and only if every or
any $f_{\lambda}$ does. Thus if some $f$ in $U$ has a zero on
$L_{\theta}$, then there exists $f_0$ in $U$ that vanishes at
$e^{i\theta}$, that is, the point in
$L_{\theta}$ on the unit circle. Thus it would suffice to show that
$f(e^{i\theta}))
\neq 0$ for all $f$ in $U$. But we can do a bit better.  

We claim that
all the zeros of $V(\bb)$ also lie in
$W_{\theta}$.  Suppose not; there exists $f = \sum c_j z^j$ in
$V(\bb)$ with a zero, $z_0$, not in $W_{\theta}$. There exist $f_n$ whose
first two coefficients are $c_0$ and $c_1$ respectively, with $f_n \to f$
coordinatewise. Since $W_{\theta}$ is regular, either $z_0$ is in
$\partial W$ or in  the complement of $W \cup \partial W$; but the latter
is impossible from $f_n \to f$.

Now let  $\bb = (\beta,\beta, \beta, \dots, \beta)$ for some $\beta > 1$, and let $U_{\beta}$ ($V_{\beta}$)
denote
$U(\bb)$ ($V(\bb)$, respectively). 

\comment
We do some special cases first. If $\theta = \pi/2$ (so $W$ is the open left half-plane), choose $f$ in $V(\beta)$
(optimal $\beta >1$ to be determined), say $f = \sum c_j z^j$. Pick a mode $k$, meaning that  $c_k \geq c_j$ for all
$j$.  Consider $i^{-k+1}f$, and take its imaginary
part: this is $\dots - c_{k-2} + c_k - c_{k+2} + c_{k+4} - \dots$, an alternating series in which the absolute value of
the terms in increasing up to $c_k$ and decreasing beyond $c_k$. Note that since $\beta > 1$, the sequence $c_0, c_1,
\dots, c_{k-1}$ is {\it strictly\/} increasing  and $c_{k+1}, c_{k+2}, c_{k+3}, \dots$ is strictly decreasing until it
hits zero. It suffices to show that
$c_k
\geq c_{k+2} + c_{k-2}$ and at least one of $c_{k\pm 4} >0$ ($N \geq 6$ will guarantee the latter). Set (if $k \geq 1$)
$r = c_{k+1}/c_k \leq 1$ and
$s = c_{k-1}/c_k \leq 1$, so that $rs = 1/\beta_1 \leq 1/\beta$. Now  assume that $N-2 \geq k \geq 2$. Then $c_{k+2}/c_k
\leq r^2/\beta$ and $c_{k-2} = s^2/\beta$. Now $\max r^2 + s^2$, subject to $0 \leq r,s \leq 1$ and $rs \leq 1/\beta$,
equals $1+ 1/\beta^2$ (attained  only when $r = 1$ and $s= 1/\beta$ or vice versa). Hence will obtain $c_k \geq c_{k+1} +
c_{k-1}$ if $1 \geq (1+ \beta^2)/\beta$, that is $\beta^3 - \beta^2 -1 \geq 0$ ($\beta \geq \rho$, well known algebraic
integer, $\rho$ is close to $1.5$). 

If $k=1$, the sequence $(c_1, c_3, c_5, \dots)$ is monotone decreasing, so it suffices that $ c_5 > 0$ (this requires
$N\geq 5$), and the same holds if $k = 0$ (this requires only that $c_4 > 0$. 

Hence we have shown that if $f= \sum c_j z^j$ with $\beta_j \geq \rho$ (where $\rho$ is the root of $x^3 -x^2 -1 = 0$
that exceeds $1$) for all
$j = 1, 2,\dots, N-1$ and
$ N
\geq 6$, then all zeros of $f$ lie in the open left half-plane.
\endcomment

Let $f = \sum_{j=0}^N c_j z^j$, with $N \in \brcs{5,6,7,\dots} \cup
\brcs{\infty}$ be an entire function with strictly positive coefficients
\st for all $1 \leq j < N$, the numbers $\beta_j:= c_j^2/c_{j+1} c_{j-1}$
are all at least as large as $\beta > 1$. We determine conditions on
$\beta$ (and to a lesser extent, on $N$), to guarantee that $f$ does not
vanish at $e^{i\theta}$, for suitable values of $\theta$.

First assume that $ N < \infty$, so we are dealing with polynomials.

We   note that if $f$ is replaced by its opposite (obtained by
reversing the order of the coefficients), the set of $\beta_j$ does not
change (only the indexing), and any zeros on the unit circle that are
zeros of the opposite function are zeros of the original. Hence in trying
to show that $f(e^{i\theta})\neq 0$ for all members of $U_{\beta}$,  we
can assume that if the mode appears at $k$ (that is, $c_k \geq c_j$ for
all $j$; since $\beta >1$, there is at most one other mode, which must be
adjacent), then $c_{k+1} \leq c_{k-1}$ (or vice versa, whichever works out
better).

For a polynomial $f$, let $Z(f)$ denote its set of zeros (multiplicities
are irrelevant for this discussion). For $0 < \theta< \pi$, let
$W_{\theta}$ be $\Set{z \in \C}{|\arg z- \pi| < \pi -\theta}$, and denote
$\pi-\theta$ by $\psi$. Define 
$$\beta_0 (\theta) = \inf \Set{\beta} {f \in
U_{\beta} \text{ implies } Z(f) \subset W_{\theta}}.$$
 (If $f \in U_4$ and
$N \geq 5$, then $Z(f)$ consists of negative real numbers, hence belongs
to $W_{\theta}$; hence $\beta_0 (\theta) \leq 4$.)

We can easily obtain lower bounds for $\beta_{0}(\theta)$ when $\pi/2 \leq
\theta < \pi$; by more difficult methods, we show these are sharp  for
each $N \geq 6$.
The functions $\RRR$ and $\RRRRR$ are defined in the introduction.
Finally, we can state the main result of this section. For $N \geq 6$,
$$
\beta_0 (\theta) = \cases 4\cos^2 \theta & \text{if $4\pi/5 \leq \theta <
\pi$}\\ 1 - 2\cos \theta & \text{if $\theta_0 \leq \theta \leq4\pi/5$}\\
\RRR(\theta) & \text{if $\theta_1 \leq \theta \leq \theta_0$}\\
\RRRRR(\theta) & \text{if $\pi/2\leq \theta \leq \theta_1$}\\
\endcases
$$ where $\theta_0 = .64\dots \cdot\pi$ (almost $2\pi/3$) is the solution to $1- 2\cos \theta = \RRR(\theta)$ and 
$\theta_1 = .53\dots \cdot\pi$ (just above $\pi/2$) is the solution to $\RRR(\theta) = \RRRRR(\theta)$. These
results extend to entire functions.

That means if  $ p = \sum_{i=0}^N c_i z^i$ ($ N \in \brcs{6, 7,
8, \dots} \cup \brcs{\infty}$) is an entire function with only nonzero
coefficients and satisfying $c_i^2/c_{i+1}c_{i-1} \geq \beta_0 (\theta)$ for some
$\theta$ but all $1 \leq i < N$, then for all zeros $z$ of $p$, $|\arg z| > \theta$, and $\beta_0 (\theta)$ 
is the smallest number with this property. For example, if $\theta = 3\pi/4$, $2\pi/3$, $\pi/2$, the respective
values of $\beta_0 (\theta)$ are $1+ \sqrt 2$, $2$, and the real root of $X^3 - X^2 -1 =0$
($((116+12\cdot93^{1/2})^{1/3} +4(116+12\cdot 93^{1/2})^{-1/3}+2)/6$; approximately $1.46557\dots$). 

Now we have the relatively easy necessary conditions.

Recall that  $\RRRRR(\theta)$ is the largest positive  root of $X^3 -
(1+a) X^2 -1 + a + a^2$ where $a = 2 \cos \psi$, when it exists, and $1$
otherwise. 

\Lem Lemma \onefiv. With $N \geq 5$ and $\pi/2 \leq \theta < \pi$ and
$\psi =
\pi - \theta$, we have that
$$
\beta_0(\theta) \geq \max\brcs{4 \cos^2 \theta, 1 + 2 \cos \psi, \RRR(\theta), \RRRRR(\theta)}.
$$

\Pf We first obtain polynomials of degrees $2$ through $5$ yielding the
lower bounds, and then show how they can be enlarged to sequences of
polynomials of degree $N$ to yield the lower bounds in all cases. For
convenience, let $a = 2\cos \psi$. Set 
$g:= z^2 + az  + 1$; its roots are 
$e^{\pm i \theta}$, and its lone $\beta_1$ is $4\cos^2 \psi$. Now set $h = g\cdot (1 + z) = z^3 + (1+a)z^2 + (1+a)z +1$,
whose $\beta$ values are $1+ a$, and obviously with roots $\brcs{e^{\pm i \theta},-1}$. 

Next, set $j = g \cdot (1 + b z
+ z^2)$ where $b$ is to be determined. This expands as $z^4 + (a + b)z^3 + (2 + ab)z^2 + (a+b)z + 1$, whose
$\beta$-values are $\brcs{(a+b)^2/(2+ab), (2+ab)^2/(a+b)^2}$. Let $b$ be a positive root of the equation $(a+b)^2/(2+ab),
=(2+ab)^2/(a+b)^2$, i.e., $(a+b)^4 = (2+ab)^3$ (if none exist, then this will correspond to the value $1$ for
$\RRR(\theta)$, so that the sole $\beta$-value of $j$ is $(a+b)^{2/3}:= \beta$. This  yields $b = \beta^{3/2} - a$,
and substituting this into the equation, we obtain $\beta^6 = (2+ a(\beta^{3/2}-a)^3$, or in other words, $\beta^2 - a
\beta^{3/2} + a^2 -2 = 0$. By going in reverse, we reconstruct $b$ from the positive real root of this quartic. Hence the
$\beta$ value of
$j$ is
$\RRR(\theta)$.

Finally, set $ k = g\cdot (1+ bz + bz^2 +z^3)$, again with $b$ to be determined. Here $k = z^5 + (a+b)z^4 + (1+b +
ab)z^3 + (1+b + ab)z^2 + (a+b)z + 1$, with at most two distinct $\beta$ values, $(a+b)^2/(1+b + ab)$ and
$(1+b+ab)/(a+b)$;
equating them as  in the previous case, we obtain the equation $(a+b)^3 = (1+ab +b)^2$; if this has a positive real root
$b$, the $\beta$ value of the corresponding choice of $k$ is $\beta:= (a+b)^{1/2}$. Then with the substitution $b =
\beta^2 -a$ yields the  equation (for $\beta$) $\beta^6 = (1+(\beta^2 -a)(1+a))^2$, and since both sides are positive,
$\beta^3 = 1 + \beta^2 (1+a) -a -a^2$, or in other words, $\beta^3 -(1+a)\beta^2 -1 + a + a^2 = 0$. Now start with this
and define $b$ to construct $k$. 

In each of the four cases we have found polynomials (of degrees two through five)
whose $\beta$ values are as indicated, and have $e^{i \theta}$ as a root. If $N \geq 5$ (or $N>5$ and
we are dealing $\RRRRR(\theta)$), let $l$ be one of $\brcs{2,3,4,5}$ and define for each integer $n$, $f_n = 1 +
z^{1}/n^2 + z^{2}/n^6 + \dots + z^{N-l}/n^{N(N-l)}$ (so that the $\beta$
values are all $1/n^2$, except if $N-l =1$). If $F$ is any polynomial
whose minimal $\beta$ value larger than $1$, it is easy to see that the
minimal beta values of the elements of the sequence 
$F\cdot f_n$ converge to that of $F$. For $F \in \brcs{g,h,j,k}$, $e^{i\theta} $ is a root of $F$ hence of $F\cdot f_n$
for all $n$, and since the degree of $F\cdot f_n$ is $N$, it follows that $\beta_0 (\theta)$ is at least the
$\beta$-value of each $F \cdot f_n$, hence is at least the limit of their $\beta$ values. 
\qed

We will frequently work with $\psi = \pi - \theta$. The proof that
$\beta_0(\theta)$ is bounded above by the right side involves overlapping
intervals on which we work with only two of the terms to be maximized at a time
(e.g., on $(3\pi/4, \pi)$, we show $\beta_0 (\theta) \leq \max \brcs{4
\cos^2 \theta, 1 - 2\cos \theta}$). The reverse inequalities (which shows
that the right side is always sharp) are obtained from tricky
multiplications, which are easy to implement, but were difficult to find.

We note the following elementary inequalities. Suppose $f = \sum c_j z^j$
belongs to $U_{\beta}$, and $k$ is a mode for the sequence $(c_j)$. Then
with $l \geq 0$ and $ j$ positive and large enough,
$$\eqalign{
\frac{c_{k+l +j}}{c_{k+l+j-1}} & = \frac{c_{k+l +j-1}}{c_{k+l+j-2}}\cdot
\frac{1}{\beta_{k+l+j-1}} \leq \frac{c_{k+l +j-1}}{c_{k+l+j-2}}\cdot
\frac{1}{\beta}  \leq \dots \cr
& \leq \frac{c_{k+l }}{c_{k+l-1}} \cdot \frac{1}{\beta^{j}} \qquad\text{and} \cr
\frac{c_{k+l +j}}{c_{k+l+j-1}} & \leq \frac{c_{k+1}}{c_{k}} \cdot \frac{1}{\beta^{j+l-1}}  \cr
}$$ and
$$\eqalign{
\frac{c_{k+l +j}}{c_{k+l}}  & = \frac{c_{k+l +j}}{c_{k+l+j-1}}
\frac{c_{k+l +j-1}}{c_{k+l+j-2}} \cdot \dots \cdot \frac{c_{k+l
+1}}{c_{k+l}} \cr
& \leq  \(\frac{c_{k+l +1}}{c_{k+l}}\)^{j}\cdot \frac1{\beta^{(j^2 +
j)/2}}\cr
}$$

Throughout, $w = e^{i\theta}$ and $\psi = \pi - \theta$.

\Lem Theorem \onesix. With fixed $N \geq 6$, we have the following.
$$
\beta_0 (\theta) = \cases 4\cos^2 \theta & \text{if $4\pi/5 \leq \theta <
\pi$}\\ 1 - 2\cos \theta & \text{if $\theta_0 \leq \theta \leq4\pi/5$}\\
\RRR(\theta) & \text{if $\theta_1 \leq \theta \leq \theta_0$}\\
\RRRRR(\theta) & \text{if $\pi/2\leq \theta \leq \theta_1$}\\
\endcases
$$
 where $\theta_0 = .64\dots \cdot \pi$  is the solution to $1- 2\cos \theta = \RRR(\theta)$ and 
$\theta_1 = .53\dots\cdot \pi$ is the solution to $\RRR(\theta) = \RRRRR(\theta)$. 
{\par}
In particular, if $f = \sum_{j = 0}^N c_j z^j$ with $N \in \brcs{6,7, 8, \dots} \cup \brcs{\infty}$
is entire with all coefficients nonzero and nonnegative,
and $\inf_{1 \leq j < N} c_j^2/c_{j+1} c_{j-1}  \geq \max\brcs{4\cos^2 \theta,1 - 2\cos \theta,\RRR(\theta),
\RRRRR(\theta)}$ for some $\theta$ in $[\pi/2,\pi)$, then all zeros of $f$ lie in the open sector $|\arg z| >
\theta$. 

\Pf We fix $6 \leq N <\infty$. By \onefou, we need only obtain conditions excluding
zeros at various $e^{i\theta}$---that is, we show that if $\beta$ is at least as large as the
right side (usually treated two at a time), and $f $ is in $U_{\beta}$, then $f(e^{i\theta}) \neq 0$. This 
yields that $\beta_0 (\theta) $ is bounded above by the right side, but we already have the reverse inequality
in \onefiv. We deal with three overlapping intervals. 

\noindent (a) $3\pi/4 \leq \theta < \pi$. Here we show that
$\beta_0 (\theta) \geq \max\brcs{4\cos^2 \phi, 1- 2\cos \phi}$.

\noindent Suppose
there exists
$f$ in
$U_{\beta}$ ($\beta$ as yet unspecified---we wish to derive conditions
that guarantee a contradiction, which will typically bound $\beta$) \st
$f(e^{i\theta}) = 0$. By replacing $f$ by its opposite if necessary, we
may assume that the coefficient to the right of the mode  is less than or
equal to the coefficient to the left of the mode. Let $k$ be the mode, so
that $c_{k+1} \leq c_{k-1}$, and consider  $-\Im
(e^{-(k-2)i\theta}f(e^{i\theta}))$. This expands (replacing $\theta$ by
$\pi -\psi$, which makes the manipulations clearer; here $0 < \psi \leq
\pi/4$) as
$$\eqalign{
[ \dots  -c_{k-4}\sin 2\psi &+c_{k-3} \sin \psi + 0\cdot c_{k-2}] \cr & +[
-c_{k-1}
\sin\psi      +c_k \sin 2\psi - c_{k+1}\sin 3\psi] \cr &+ [c_{k+2} \sin
4\psi - c_{k+3 }\sin 5\psi + \dots];\cr
}$$
we interpret $c_{\text{negative}} = 0$. We will analyze this in three
parts, the middle three terms, $ -c_{k-1} \sin\psi      +c_k \sin 2\psi -
c_{k+1}\sin 3\psi$, the right tail, $c_{k+2} \sin 4\psi - c_{k+3 }\sin
5\psi + \dots$, and the left tail, $ \dots  -c_{k-4}\sin 2\psi +c_{k-3}
\sin \psi $. We will show that if $\beta \geq \max\brcs{4\cos^2 \psi, 1 + 2\cos
\psi}$, then all three are nonnegative, and at least one of the tails is
positive (of course, depending on $k$ and $N$, the left or right tail
might not even exist). The computation of the middle term brings out the
connection with the necessary conditions, while the tails just require
estimates on the rate of decay (which is very rapid).

Set $r = c_{k-1}/c_k$ and $s = c_{k+1}/c_k$, so that $0 \leq s \leq r \leq
1$ and $rs \leq 1/\beta$. To show the middle term is nonnegative, it
suffices to show that
$$
\max \Set{r\sin \psi  + s\sin 3\psi}{0 \leq r \leq s \leq 1, \ rs \leq
1/\beta} \leq \sin 2\psi.
$$
It is immediate that the only two locations  for the maximum value of the left side occur at
$(r,s) = (\beta^{-1/2},\beta^{-1/2})$ and $(1,1/\beta)$.

The value at the former leads to
$\beta^{-1/2}(\sin \psi + \sin 3\psi) \leq \sin2\psi$.  Since $\sin \psi +
\sin 3\psi = 2 \sin 2\psi \cos \psi $, we obtain $\beta^{-1} \leq 1/2\cos
\psi$, or $\beta \geq 4 \cos^2 \psi$ (that was easy).

The latter leads to $\beta^{-1} \leq (\sin 2\psi - \sin \psi)/\sin 3\psi$.
From $\sin 2\psi= 2\sin \psi \cos \psi$ and $\sin 3\psi = 3\sin \psi -
4\sin^3 \psi = \sin \psi (4\cos^2 \psi - 1)$, we see that  the inequality
is equivalent to $\beta^{-1} \leq (2\cos \psi -1)/(4\cos^2 \psi -1) =
1/(2\cos \psi +1)$.

Thus if $\beta \geq \max\brcs{4 \cos^2 \psi, 1+2\cos \psi}$  (that is, both
inequalities occur), then the middle  term is nonnegative.

Now we look at the right tail. We begin by supposing that $\pi/(K+1) \leq
\psi < \pi/K$ for some positive integer $K \geq 5$ (we have to consider
the case that $K =4$ separately). Then $\sin 4 \psi, \sin 5\psi, \dots ,
\sin (K+1)\psi$ are all nonnegative, and thus the sequence that appears in
the expansion of the tail, that is,
$$c_{k+2}\sin 4 \psi, -c_{k+3}\sin 5\psi, c_{k+4}\sin 6\psi, \dots,
(-1)^{K+1}c_{k+K-1}\sin (K+1)\psi
$$ is alternating. We will show that the
absolute values are decreasing, and so get a good lower bound  for this
portion of the tail. The remainder of the tail is so small that it is easy
to deal with.

The ratio of   absolute values of consecutive terms of $(-1)^j c_{k+j}
\sin (j+2)\psi$ ($j = 2,3, \dots , K+1$) is $(c_{k+j +1}/c_{k+j})(\sin
(j+3)\psi/\sin (j+2)\psi)$.
The left factor is bounded above by $s/\beta^j\leq \beta^{-j-1/2}$, and the
right term is bounded above by $2$ (note that for $j > (K-2)/2$, each
$\sin (j+3)\psi/\sin (j+2)\psi$ is  less than one. On the interval for
$\psi$, $\beta$ is at least $\max \brcs{4\cos^2 \pi/{K},1+ 2\cos \pi/K} >
2.5$. Since $\beta^{-j-1/2} < 1/2$, the sequence is descending.

Hence the sum of the alternating sequence $\sum_{j=2}^{K-1}
(-1)^jc_{k+j}\sin (j+2)\psi$ is bounded above by $c_{k+2} \sin 4\psi -
c_{k+3}   \sin 5\psi$. This can be rewritten as $c_{k+2} \sin 4\psi (1-
(c_{k+3}/c_{k+2})(\sin 5\psi/\sin 4\psi))
\leq c_{k+2} \sin 4\psi  (1- (\sin 5\psi/\sin 4\psi) \beta^{-5/2}$. The
remainder of the tail  can be bounded in absolute value by
$$\eqalign{
\sum_{j=0}^{\infty}c_{k+ K+j}|\sin (K+j+2)\psi| & \leq
\sum_{j=0}^{\infty}c_{k+ K+j} \cr
& \leq c_{k+2} \sum_{j=0}^{\infty}\frac {c_{k+ K+j} }{c_{k+2}}  \cr
& \leq c_{k+2} \sum_{j=0}^{\infty}\beta^{-((K+j)^2 + (K+j) -6)/2}  \cr
}$$
The series is  $\beta^{(-K^2 -K +6)/2} + \beta^{(-K^2 -3K +5)/2} +
\dots $ which is bounded above by $(5/4)\beta^{(-K^2 -K +6)/2}$. So all we need
is that $\sin 4\psi (1- (\sin 5\psi/\sin 4\psi)\beta^{-5/2}) >
(5/4)\beta^{(-K^2 -K +6)/2}$ for $K \geq 5$ and $\psi $ in the interval.
Sufficient is that $\sin 4\psi > 2.4^{(-K^2 -K +6)/2}$. Since $\sin 4\psi >
\psi > \pi/K$ on this interval, the result follows easily.

Now suppose (still dealing with the right tail) that $\pi/5 \leq \psi \leq \pi/4$. In this case, $\sin 5 \psi$
is negative, and the minimum of $c_{k+2} \sin 4\psi - c_{k+3} \sin 5\psi$
occurs at the right endpoint (easy to check, since $c_{k+3}/c_{k+2} \leq
s/\beta^2 < 1/2$), which is $\sqrt {3} c_{k+3}/2$. The remainder of the
tail is bounded by $\sum_{j\geq 4} c_{k+j}$, which is dealt with as above.

The treatment of the left tail is similar, but a little simpler. Again assume that $\pi/(K+1) \leq \psi \leq
\pi/K$. The sequence $\sin \psi, -\sin 2\psi, \dots, (-1)^{K-1} \sin K\psi$  is alternating and $c_{k-j-2}\sin j
\psi$ (for $j = 1,\dots, K$) is descending, so 
$$\eqalign{
\sum_{j=1}^K (-1)^{j-1}c_{k-j-2} \sin j\psi & \geq c_{k-3} \sin \psi - c_{k-4} \sin 2\psi\cr
& \geq c_{k-3} (\sin \psi -  \sin 2\psi/\beta^3) = c_{k-3}\sin \psi (1 - 2\beta^{-3}\cos \psi) \cr
& \geq c_{k-3}\sin \frac{\pi}{K+1} (1-2\beta^{-3}). \cr
}$$ 
The remaining terms in the left tail are bounded above by $\sum_{j=0}c_{k-K-2-j}$, which as before, is bounded
above by $c_{k-3} \( \beta^{(-K^2 + K)/2} + \beta^{(-K^2 - K)/2} + \dots \)$. All that remains is to verify that 
$(1-2\beta^{-3}) \sin \pi/(K+1) > 1.1\beta^{(-K^2 + K)/2}$ for $\beta > 2.4$ and $K \geq 4$, which is easy.

\noindent (b) $3\pi/5 \leq \theta \leq 3\pi/4$, i.e., $\pi/4
\leq
\psi \geq 2\pi/5$.  Here we show that
$\beta_0(\theta)
\leq
\max\brcs{1+ 2\cos \psi, \RRR(\theta)}$ 

\noindent Recall that
 $\RRR(\theta)$ is the largest positive real zero  of  (what is
effectively) a quartic, $X^{2} + X^{3/2} 2 \cos \theta+ 2\cos 2\theta
=0$; if none exist, define $\RRR(\theta) = 1$.

On this interval, $\sin \psi$ and $\sin 2\psi$ exceed zero, and $\sin 4\psi$ is nonpositive; $\sin 5\psi < 0$ if $\psi
< 2\pi/5$. Let
$k$ be the mode of
$\brcs{c_k}$ and assume that $c_{k-1} \leq c_{k+1}$. Consider the middle terms, $-c_{k-1}\sin \psi + c_k \sin 2\psi -
c_{k+1}
\sin 3\psi + c_{k+2}\sin 4\psi$  of
$\Im w^{-k+2} f(w)$; set
$r = c_{k-1}/c_{k}$ and $s = c_{k+1}/c_k$, so that $0 \leq r \leq s \leq 1$ and $rs \leq 1/\beta$. We derive conditions
(on $\beta$) to guarantee that $-r \sin \psi + \sin 2\psi -  s\sin 3\psi  + 4 s^2/\beta \sin 4\psi$ is nonnegative
(which is sufficient for nonnegativity of the middle term  since $c_{k+2} \leq s^2 c_{k}/\beta$). This amounts to
$r + s (4\cos^2\psi -1) + s^2 |8\cos^3 \psi - 4\cos\psi|/\beta \leq 2\cos \psi$. Now the maximum value of the left side
occurs at either
$(r,s) =
\beta^{-1/2}(1,1)$ or at $(1/\beta,1)$. The former yields a maximum value of $\beta^{-1/2}(4\cos^2\psi) -
4\cos \psi (2\cos^2\psi -1)\beta^{-2}$, which is less than or equal to the right side if and only if $\beta^2
-\beta^{3/2} 2\cos \psi + 2\cos 2\psi \geq 0$, that is, $\beta \geq \RRR(\theta)$ (in converting between $\theta$ and
$\psi$, the middle term is multiplied by $-1$, but not the constant term).

At      $(r,s) = (1/\beta,1)$, we require
$\beta^{-1} (1 - 4\cos \psi(2\cos^2 \psi -1)) \leq 1 + 2\cos \psi- 4\cos^2\psi$, that is $\beta^{-1} \leq (1 + 2\cos
\psi- 4\cos^2\psi)/(1 - 4\cos \psi(2\cos^2 \psi -1)) = 1/(1+2\cos \psi)$. Hence $\beta \geq \max \brcs{1+2\cos
\psi,
\RRR(\theta)}$ is sufficient to guarantee that this middle cluster of terms is nonnegative. Now we deal with the tails.

Now we deal with the right tail. The tail begins $-c_{k+3}
\sin 5\psi + c_{k+4} \sin 6 \psi - 7 c_{k+5} \sin 7 \psi + \dots$. The first subcase is $\pi/4 \leq \psi \leq
\pi/3$.  On this interval, the leading term,
$-c_{k+3}
\sin 5\psi$ is at least $c_{k+3}/\sqrt 2$. The rest of the tail is bounded in absolute value rather crudely by
$\sum_{j\geq 1} c_{k+3+j}$. Obviously 
$$\eqalign{
\sum_{j\geq 1} c_{k+3+j} & = c_{k+3} \sum_{j\geq 1} \frac{c_{k+3+j}}{c_{k+3}} \leq c_{k+3} \sum_{j\geq 1}
\frac{s^j}{\beta^{(j^2 + 5j)/2}}\cr
& = c_{k+3} \( \frac 1{\beta^{3}} + \frac1{\beta^7} +\frac1{\beta^{12}} + \dots  \).\cr 
}$$
With $\beta \geq 2^{1/2}$, we deduce $-c_{k+3} \sin5\psi > \sum c_{k+3 +j}$, hence  the right tail is positive.

The next subcase assumes $\pi/3 \leq \psi \leq 2\pi/5$. Then the leading two terms are nonnegative, and 
$$\eqalign{
-c_{k+3}
\sin 5\psi + c_{k+4} \sin 6 \psi &\geq c_{k+4} (\sin 6\psi - \sin 5\psi) > .9 c_{k+4}\cr
}$$
Now the rest of the tail is bounded in absolute value by $\sum_{j\geq 1} c_{k+4+j}$, and the same technique as in the
previous subcase yields that the right tail is positive.

For the left tail, $c_{k-3} \sin \psi - c_{k-4} \sin 2\psi + c_{k-5} \sin 3\psi - \dots$, the leading term is at least
$c_{k-3}/\sqrt 2$, and the same argument as in the first subcase of the right tail (but without requiring further
restrictions on $\psi$) will work. 

\noindent (c) $13\pi/36 \leq \psi \leq \pi/2$. 

\noindent Here we show that $\beta_0
(\theta) \leq \max \brcs{\RRR(\theta), \RRRRR(\theta)}$.

\noindent We note the following identities
$$\eqalign{
\sin \frac{3\psi}2 & =\sin \frac{\psi}2 \(1 + 2 \cos \psi \)\cr
\sin \frac{5\psi}2 & =\sin \frac{\psi}2 \(4\cos^2 \psi+ 2 \cos \psi -1
\)\cr }$$

We may assume the mode appears at $k$ and $c_{k-1} \leq c_{k+1}$.
Consider $-\Im w^{-k+3}f(w)$, which expands as
$$\eqalign{
[\dots -c_{k-4}\sin 7\psi &+c_{k-3}\sin 6\psi ]\cr &+ [- c_{k-2}\sin 5\psi
+ c_{k-1}\sin 4\psi -c_{k} \sin 3\psi + c_{k+1}\sin 2\psi  - c_{k+2}\sin
\psi ]\cr 
& + [c_{k+4}\sin \psi - c_{k+5}\sin 2\psi + \dots] \cr
}$$
(Note that   $\sin 4\psi$ and $\sin 3\psi$ are both nonpositive and $\sin
6\psi$ is positive  on $\pi/3 < \psi \leq
\pi/2$.) As usual, begin with the middle thing; with $r = c_{k-1}/c_{k}$ and $s = c_{k+1}/c_{k}$, it is sufficient to
determine conditions on $\beta$ to guarantee that 
$$
- \frac{r^2}{\beta}\sin 5\psi + \frac r{\beta}\sin 4\psi - \sin 3\psi + \frac s{\beta}\sin
2\psi  - \frac{s^2}{\beta}\sin \psi \geq 0,
$$
subject to the constraints $0 \leq r \leq s \leq 1$ and $rs \leq 1/\beta$. It is again easy to check that the minimum
value occurs at the the vertices of the domain, $(r,s) = \beta^{-1/2}(1,1)$ and $(\beta^{-1},1)$. Evaluating at the first
yields
$$
\frac{-1}{\beta^2} (\sin 5\psi + \sin \psi) + \frac{1}{\sqrt{\beta}} (\sin 4\psi + \sin \psi) - \sin 3\psi.
$$
Noting that $-\sin 3\psi > 0$ (if $\pi/3 < \psi \leq \pi/2$), we may divide by $(-\sin 3\psi)\beta^2$, and expanding the
sums of sines, we obtain the equivalent inequality,
$\beta^2 - 2\beta^{3/2}\cos \psi + 2\cos 2\psi \geq 0$, for which $\beta \geq \RRR(\theta)$ is sufficient.

At the point $(r,s) = (1/\beta,1)$, we obtain $\sin 2\psi - \sin 3\psi + \beta^{-1}\cdot(\sin 4\psi -\sin \psi) -
\beta^{-3}\cdot \sin 5\psi$. Converting the differences of sines and noting that $-\cos 5\psi/2 > 0$, we can divide this
by $-2\cos (5\psi/2) \sin \psi/2 \beta^{-3}$, and obtain the equivalent, 
$$
\beta^3 - \frac{\sin 3\psi/2}{\sin \psi/2}\beta^2 + \frac{\sin 5\psi/2}{\sin \psi/2} \geq 0.
$$
For this, $\beta \geq \RRRRR(\theta)$ is sufficient, by the identities above.

Hence if $\beta \geq \max\brcs{\RRR(\theta),\RRRRR(\theta)}$, the middle part is nonnegative; this was under the
assumption that $\pi/3 < \psi \leq \pi/2$. Now we deal with the tails. 

For the left tail, we first consider $13\pi/36 \leq \psi \leq 3\pi/7$. Then $13\pi/6 \leq 6 \psi \leq 18/7\pi$,
so $\sin 6 \psi \geq 1/2$. The rest of the tail is bounded in absolute value  by $\sum_{j\geq 0}
c_{k-4-j}$, which as usual is  bounded above by $c_{k-3} ( \beta^{-3} + \beta^{-6} + \beta^{-10} + \dots$. Since
$\beta > \sqrt 2$, the factor is less than $1/2$, and so the left tail is positive.

Next, assume $3\pi/7 \leq \psi \leq \pi/2$. Then $3\pi \leq 7 \psi \leq 7\pi/2$, so $ -\sin 7\psi \geq 0$. Hence
$$
c_{k-3} \sin 6 \psi - c_{k-4} \sin 7 \psi \geq c_{k-4} (\sin 6 \psi -  \sin 7 \psi) = 2c_{k-4} \sin \frac{\psi}2
\left|\cos \frac{13\psi}2\right| > .85 c_{k-4}.
$$
The rest of the tail is bounded above in absolute value by $\sum_{j\geq 0} c_{k-5 -j} \leq c_{k-4} (\beta^{-4} +
\beta^{-7} + \beta^{-11} + \dots)$, which is a lot less than $.85 c_{k-3}$ (when $\beta > \sqrt 2$). 

For the right tail, we note that $\sin \psi > .95$ on the entire interval $[3\pi/7, \pi/2]$ ($\sin 3\pi/7 =
.9749\dots$), and the rest of the tail is bounded above by $c_{k+4} (\beta^{-4} + \beta^{-7} + \dots) < .5
c_{k+4}$.

Hence we have that (with fixed finite $N \geq 6$), $\beta_0 (\theta)$ is given by the right side. When $N$
 is finite, the final statement is a restatement of the definition of $\beta_0 (\theta)$. When $N $ is
infinite, we note that the finite truncations converge uniformly on compact sets to $f$, and each of these
truncations satisfy the same conditions, hence their zeros lie within $W_{\theta}$, and thus so do all the zeros
of $f$.
\qed

Owing to the awkward definitions of $\RRR(\theta)$ and $\RRRRR(\theta)$, it is worthwhile discussing the
function inverse to $\Arrow \beta_0;[\pi/2,\pi).[1.46557\dots,4]$ ($\beta_0$ is of course defined on all of
$(0,\pi)$, but we only have an exact formula available on $([\pi/2,\pi))$. 

Let $P_N^{++}$ denote the set of polynomials of degree $N$ all of whose coefficients are strictly positive;
obviously $U_{\beta} \subset P_N^{++}$ for all $\beta > 1$. Taking, as usual the branch of $\arg z$ given by
$-\pi < \arg z \leq \pi$, define 
$$
\Arrow T;P_N^{++}.(0,\pi] \qquad\text{by} \qquad T(f) = \inf\Set{|\arg z|}{z \in Z(f)}.
$$
Since polynomials in $P_N^{++}$ have no positive real numbers as zeros, $T$ is well-defined. A simple
RouchŽ-convergence argument shows that $T$ is continuous. 

Define $\Arrow \Theta;(1,4].[0,\pi]$ via 
$$
\Theta (\beta) = \sup \Set{\theta \in (0,\pi]}{Z(f) \subset W_{\theta} \text{ for all $f\in U_{\beta}$}}. 
$$
It is easy to check that $\Theta = \beta_0^{-1}$, and moreover, $\Theta (\beta) = \inf \Set{T(f)}{f \in
U_{\beta}}$. On the interval $[1.46557\dots,4]$ (the left endpoint is the real root of $X^3 - X^2 - 1 =0$), we
have (where $\gamma = (1+ \sqrt 5)/2$ is the golden ratio)
$$
2 \cos \Theta(\beta) = \cases -\sqrt{\beta}& \text{if $\gamma^2 \leq \beta \leq 4$}\\ 
1-\beta & \text{if $1.57762\dots \leq \beta \leq \gamma^2$} \\
-\beta^{3/2} + \sqrt{\beta^3 - 4\beta^2 + 8} & \text{if $1.52334\dots \leq \beta \leq 1.57762\dots$} \\ 
1-\beta^2 + \sqrt{(1+ \beta^2)^2 + 4(1-\beta^3)} & \text{if $1.46557\dots \leq \beta \leq 1.52334\dots$} \\
\endcases
$$
\vskip 6pt
\noindent {\it If $\theta < \pi/2$.} When $\theta < \pi/2$, especially as $\theta \to 0$, new phenomena occur.
The first is that the dependence  on $N$ (which was barely noticeable up to this point) becomes more marked. It is
worthwhile giving an equivalent form of  [H, Corollary 1.3].

\Lem Proposition \onesev. Suppose that $q$ is a monic polynomial of degree $n-1$ or less, and has real nonnegative
coefficients which  form a unimodal sequence. If $q$ has a zero at the single point $\exp(2\pi i/n)$,
then  $g = 1 + z + z^2 + \dots + z^{n-1}$. 

In this result, $q$ is assumed to vanish only at the single primitive root of unity,
$e^{2\pi i/n}$, from which (together with unimodality) we deduce that it vanishes at {\it all\/} $n$th
roots of unity. From this, in order to obtain meaningful results about $\beta_0 (\theta)$ for $\theta = 2\pi/n$
(and for values close to this), we must assume that the corresponding $N$ is at least $n$, and likely at least
$3n/2$. For that reason, we redefine $\beta_0 (\theta)$ in what follows to be the $\liminf_{N \to \infty}$ of the
previously-defined $\beta_0 (\theta)$. It is likely that for fixed $\theta$, the sequence is
ultimately stationary. 

Lower bounds for the values of $\beta_0$ (for $\theta < \pi/2$) seem intractible at the moment---they likely
involve multiplication by polynomials of increasing degree (see the proof of \onefou). On the other hand, as
$\theta \to 0$, we can obtain asymptotic estimates for $\beta_0 (\theta) -1$; specifically, $\beta_0 (\theta) - 1
\leq 16 \theta^2\ln 2/\pi^2$ (the constant, $16\ln 2/\pi^2$, is rather flabby, and doubtless can be improved), and
it is relatively easy to see that
$\beta_0 (\theta) - 1
\geq
\theta^2/4\pi^2$. Hence
$\beta_0 (\theta) -1 $ is bounded above and below by a multiple of $\theta^2$.

As $\beta_0 $ is monotone, modulo a bit of fiddling with the scalar multiple, to prove the first statement,
$\beta_0 (\theta) -1 \leq K \theta^2$, it suffices to do it with $\theta = \pi/n$, where $n$ is a positive
integer. We illustrate the case that $4$ divides $n$; the other cases are very similar.

\noindent {\it $ \theta = \pi/n$ and $4$ divides $n$.}
Assume $\beta^{n^2/2} \geq 2$ (this will be subsumed by a stronger condition), mode at $k$; form $\Im
w^{-k+n/2}f(w)$; then the middle clump in the expansion consists of about $3n$ terms,
$$
\sum_{j = -n/2}^{n/2} c_{k+ j} \sin (n/2 +j)\pi/n - \sum_{j=-3n/2}^{-n/2
-1} c_{k+j} \sin (3n/2 + j)\pi/n - \sum_{j=n/2+1}^{3n/2 -1} c_{k+j} \sin
(3n/2 - j)\pi/n
$$
where we have indexed the arguments of the sine in order to ensure that
the sine terms are nonnegative. Now we show nonnegativity of each of the
following terms,
$$\eqalign{
c_{k-n/2 +1} \sin \pi/n &- (c_{k- 3n/2 +1} \sin \pi/n + c_{k-3n/2 +2} \sin
2\pi/n)\cr
c_{k-n/2 + 2} \sin 2\pi/n &- (c_{k- 3n/2 +3} \sin 3\pi/n + c_{k-3n/2 +4}
\sin 4\pi/n)\cr
\dots & \dots \cr
c_{k-n/2 + l} \sin l\pi/n & - (c_{k- 3n/2 +2l-1} \sin (2l-1)\pi/n +
c_{k-3n/2 +2l} \sin 2l\pi/n)\cr
\dots & \dots\cr
c_{k-n/4} \sin \pi/4 & - (c_{k- n -1} \sin (n/2-1)\pi/n + c_{k-n } \sin
\pi/2)\cr
}$$
We note that 
$$\eqalign{
\frac{c_{k -3n/2 + 2l-1}}{c_{k-n/2 + l}} &\leq \frac1{\beta^{(n-2l
+1)(n-2l+2)/2}}\cr
\frac{c_{k -3n/2 + 2l}}{c_{k-n/2 + l}} &\leq \frac 1{\beta^{(n-2l
+2)(n-2l+3)/2}}\cr 
}$$
 Hence if $e(l)= n-2l)^2 + 3(n-2l) + 2)/2$, we need only
show that for $l = 1,\dots, n/4$,
$$
\beta^{e(l)} \geq \frac{\sin(2l-1)\pi/n + \sin 2l\pi/n}{\sin l\pi/n}.
$$
The right side is always less than $4$, so sufficient is $\beta^{e(l)}
\geq 4$, as occurs if $\beta^{n^2/8 + 3n/4 + 1} \geq 4$. Hence $\ln \beta
\geq 2\ln 2/(n^2 /8 + 3n/4)$, or merely $\ln \beta \geq (16\ln 2)/n^2$ is
sufficient for this collection of inequalities.

The next batch of inequalities is treated similarly, but there is a slight difference. 
We consider the following,
$$
c_{k-n/4 +l} \sin \(\frac{\pi}4 + \frac {\pi l}n\) - c_{k- n +2l-1} \sin \(\frac{\pi}2 + \frac{\pi(2l-1)}n\) -
c_{k-n +2l}
\sin \(\frac{\pi}2 + \frac{2\pi l}n\),  \quad l = 1, 2, \dots, \frac n4 -1.
$$
The claim, as above, is that when $\beta$ sufficiently large, these are all nonnegative. We have that
$c_{k-n+2l-1}/c_{k-n/4 + l} \leq \beta^{-e(l)}$, where $e(l) = 15 n^2/32 - 7nl/4 + n + 3(l^2 -l)$. The smallest
value of $e(l)$ occurs (over the range $1 \leq l \leq n/4 -1$) when $l = n/4 - 1$, and we see that  $e(l ) >
n^2/8$. We also have that the ratio of the sines is bounded above by $\sqrt 2$ (some of the
flabbiness creeps in here), so sufficient for all the inequalities to hold is that $\beta^{n^2/8} > 2\sqrt 2$,
that is,
$\ln \beta > 12 \ln 2/n^2$ and in fact these are all strict. 
We have $c_k$ left over (which does not happen in the other cases, that is, when $n$ is not divisible by four).

Thus sufficient for the middle block of $3n$ or so terms to be nonnegative, it is sufficient that $\ln \beta \leq
16 \ln 2/n^2$. The remainder of the expansion is block alternating, and monotone in each of the $n$ positions,
hence the outcome is nonnegative. 

In particular, $\beta_0 (\pi/n) < \exp(16 \ln 2/n^2) \sim 1 + 16\ln 2/n^2 $. 

To give a rough upper bound for $\beta_0 (\pi /n)$, we simply note that $C_{2n}:= (1-z^{2n})/(1-z) = \sum_{0\leq j
\leq 2n-1}z^{j}$ has $\exp (2\pi i/2n)$ as a root and $C_n C_{2n}$ has sequence of coefficients
$(1,2,3,\dots,n,n,n-1,\dots, 1)$. The minimum $c_j^2/c_{j-1}c_{j+1}$ occurs when $j = n-2$, i.e.,
$(n-1)^2/n(n-2) = 1 + 1/n(n-2)$. Hence $\beta_0 (2\pi/n) -1 > 1/n^2$, so $\beta_0 (\pi/n)  -1 > 1/4n^2$. 

Both the upper and lower estimates were obtained rather sloppily, and it is very unlikely that either one is even
close to being sharp.
\vskip 5pt

\vskip 10pt\noindent {\bf Section 2. Constant quadratic ratios} \vskip 4pt
\noindent Define  for $N$ a positive integer, the polynomial of degree $N$,
$f_{b,N} = 
\sum_{j=0}^N x^{j} b^{-j(j+1)}$; if $N$ is infinite, the resulting
series is entire. This satisfies the property that 
$c_j^2/c_{j+1}c_{j-1} = b^2$ for all $1\leq j <N$. We will
determine to within $10^{-24}$, the minimum of the $b$ \st $f_{b,N}$
has only real zeros, for sufficiently large $N$, and also for
infinite $N$.
\def\Op{^{\text{op}}}

We define the {\it opposite\/} of a polynomial of degree $N$ to be the 
polynomial with coefficients written in reverse order; explcitly, $f\Op (x) 
= x^N f(x^{-1})$. A polynomial is {\it symmetric\/} (or {\it 
self-reciprocal\/}) if $f\Op = f$. If $f$ is symmetric of degree $N$, then 
$f(x) = x^Nf(x^{-1}) $, and in particular, the set of zeros of $f$ are 
closed under the operation $w \mapsto 1/w$. Every real polynomial $f = 
\sum_{j=0}^N c_j z^j$ with all $c_j^2/c_{j+1} c_{j-1}$ equal, say to $\beta 
 > 0$ for $1 \leq j \leq N-1$ can be reparameterized so as to be
symmetric.  For $N$ odd (that is, an even number of coefficients),
up to scalar multiple, the distribution of  coefficients of a
symmetric polynomial satisfying this condition
$(\dots
\beta^{-3} \ \beta^{-1} \ 1\ 1 \ \beta^{-1} \ 
\beta^{-3} \ \dots \ )$ (the exponents are triangular numbers), while for 
$N$ even, the distribution  is $(\dots  \beta^{-9/2} \ \beta^{-2} \ 
\beta^{-1/2} \     1 \  \beta^{-1/2} \ \beta^{2} \ \beta^{-9/2} \ \dots)$ 
(the exponents are half-squares).

Form ${\Cal F} \equiv {\Cal F}_{b,N} (x) = f_{b,N}(xb^{N+1})$ (for
which 
$c_j^2/c_{j+1}c_{j-1} = b^{2}$ for all relevant $j$) and write $N = 2r-1$ 
($N$ odd) or $N = 2r$ (if $N$ is even).  We see that the ratio of
the  coefficient of $x^{r}$ in $\Cal F$ to that of $x^{r-1} $ is $1$
if $N$ is  odd, and is $b$ if $N$ is even; since $\Cal F$ also has
the property that  all the ratios $c_j^2/c_{j+1} c_{j-1}$ are equal
to $b^{2}$, this enough to  guarantee that $\Cal F$ is symmetric. In
particular, for all real $x > 0$,  
$\sign {{\Cal F}_{b,N} (-x^{-1})} = (-1)^N \sign {{\Cal F}_{b,N} (-x)}$.

For $l$ a  positive integer less than or equal to $N$, consider 
$f_{b,l}$; its list of consecutive coefficients is just an initial segment 
of that of $f_{b,N}$.

\Lem Proposition \throne. Suppose that $b > \sqrt{3}$, that $l$ is a positive
integer, and that  $N \geq 2l$ is an integer. {\par}
\item{(a)} If $l$ is  even   and 
there exists a positive real number $x_1$ \st $1 < x_1 < b^4$ and
$f_{b,l}  (-x) \leq 0$, then for all $b' \geq b$, $f_{b',N}$ has all
of its roots real  and simple.{\par}
\item{(b)} If $l$ is odd and for all $x$ in $(0,b^4)$, $f_{b,l} (-x)
\geq 0$, then for all $\sqrt 3 < b' < b$, $f_{b',N}$ has nonreal
roots.

\Pf (a) We show that $f_{b,N}$ has $N+1$ sign changes along the negative reals, 
implying the result. Write $N = 2r-1$ or $2r$, depending on whether $N$ is 
odd or even. First, we will define a sequence $0 = x_0 < x_1 < x_2 < \dots < 
x_{r-1}$ of $r$ positive real numbers \st $\sign{f_{b,N}(-x_j) }= (-1)^j $, 
and then we will use the reparameterization to a symmetric polynomial, to 
show that this set can be extended to a strictly increasing sequence of 
$N+1$ positive real numbers with the same property.

To begin, we show that $f_{b,N}(-x_1) < 0$ (where $x_1> 1$ satisfies the 
conditions in the statement of this lemma). The sequence of coefficients of 
$f_{b,N}$ is strongly unimodal, hence $\brcs{x^j b^{-j(j+1)}}$ is strongly 
unimodal and thus unimodal for any choice of $x > 0$. Set $q_j = x_1^j 
b^{-j(j+1)}$. Then $q_{l+2}/q_{l+1} = x_1 b^{-2l+4}$; since $l \geq 2$, 
sufficient that this ratio be less than $1$ is that $x_1 < b^{8}$, which is 
more than satisfied. Hence the sequence $\brcs{q_j}_{j \geq l+1}$ is 
monotone descending (since the whole sequence is unimodal), and since the 
first two terms are strictly decreasing, we see $\sum_{j \geq l+1} (-1)q^j$ 
has the same sign as $(-1)^{l+1}q_{l+1}$. This is $-1$ (as $l$ is even). On 
the other hand, $f_{b,N} (-x_1) = f_{b,l}(-x_1) + \sum_{j \geq l+1} 
(-1)q^j$,
which is thus negative.

For $k = 0$, set $x_0 = 0$, so that $f_{b,N}(x_0) = 1 = (-1)^0$.
Now for each $k$ with $2 \leq k \leq N-l$, define $x_k = x_1 b^{2k-2}$. 
Since $b >1$, we have that $\brcs{x_k}_{k\geq 1}$ is strictly increasing, 
and since $N \geq 2l$, there are at least $N/2$ of them. Set $a(j) = 
j(j+1)$. Now we note the following self-replicating property (line 4
of the display below) of these polynomials. Set $x = Xb^{\alpha}$
where
$\alpha =  2k-2$
$$\eqalign{
f_{b,N}(x) & = f_{b,k-2} (x) + b^{-a(k-1)} x^{k-1} \sum_{j=0}^l x^j 
b^{a(k-1)- a(j+k-1)}
+ x^{k+l}\sum_{j=0} x^j b^{-a(j+k+l)} \cr
&=f_{b,k-2} (x) + b^{-a(k-1)} x^{k-1} \sum_{j=0}^l X^j b^{j\alpha +
a(k-1)-  a(j+k-1)}
+ x^{k+l}\sum_{j=0} x^j b^{-a(j+k+l)}  \cr
& =  f_{b,k-2} (x) + b^{-a(k-1)} x^{k-1} \sum_{j=0}^l X^j b^{-j(j+1 - \alpha 
+ 2k-2)}
+ x^{k+l}\sum_{j=0} x^j b^{-a(j+k+l)}  \cr
& =  f_{b,k-2} (x) + b^{-a(k-1)} x^{k-1} f_{b,l} (X)
+ x^{k+l}\sum_{j=0} x^j b^{-a(j+k+l)}   \cr
& =  f_{b,k-2} (x) + b^{-a(k-1)} x^{k-1} f_{b,l} (xb^{2-2k})
+ x^{k+l}\sum_{j=0} x^j b^{-a(j+k+l)}; \quad \text{since $l$ is even,} \cr
f_{b,N}(x) & = f_{b,k-2} (-x) + (-1)^{k-1}C_k f_{b,l} (-xb^{2-2k}) + 
(-1)^{k}x^{k+l}\sum_{j=0} (-x)^j b^{-a(j+k+l)}, \cr
}$$
where $C_k > 0$. Substitute $x = x_k = x_1 b^{2k-2}$. The middle term is 
then $(-1)^{k-1} C_k f_{b,l}(-x_1)$; by hypothesis, this is either zero or 
has sign $(-1)^k$. For the two tails (left and right; at least one of them 
must be nonempty), set $q_j = (x_k)^j b^{-a(j)}$. Again, since 
$\brcs{b^{-a(j)}}_{j=0}^N$ is strongly unimodal, so is the reparameterized 
sequence $\brcs{q_j}$. If $k = 2$, the left tail consists of a single term 
$1$, which has sign $(-1)^k$. If $k \geq 3$, the left tail sums to 
$\sum_{j=0}^{k-2} (-1)^jq_j$.  We note that $q_{k-2}/q_{k-3} = x_k b^{6-2k} 
= x_1 b^{4} >1$. Hence the sequence  $\brcs{q_j} _{j=0}^{k-2}$ is 
increasing, with the last difference strict. Hence $\sign{ \sum_{j=0}^{k-2} 
(-1)^j q_j}$ is the sign of the largest term, $(-1)^{k-2} q_{k-1}$, i.e., $(-1)^k$.

The right tail is treated similarly. If $k+l = N$, there is only one term, 
and its sign is   $(-1)^{k+l}= (-1)^k$; otherwise, suppose $k+ l < N$. We note that 
$q_{k+l}/q_{k+l+1} = x_k^{-1} b^{2(k+l +1)} = x_1^{-1} b^{2l } \geq b^4/x_1 
 > 1$. As an interval in a strongly unimodal sequence, it follows that 
$\brcs{q_j}_{j \geq k+l}$ is descending, and thus the sign of $\sum_{j \geq 
k+l} (-1)^j q_j$ is that of the initial term, i.e., $(-1)^{k+l} = (-1)^k$.

Hence each of the  three parts is either zero or has sign that of $(-1)^k$, and at 
least one of the three parts is not zero. Hence $\sign {f_{b,N}(-x_k)}  = 
(-1)^k$.

The $k$ for which this is valid include all $k$ with $k \leq N-l$, and since 
$N \geq 2l$ by hypothesis, this is true for all $k \leq N/2$. Now consider 
the symmetric form of $f_{b,N}$, given above as ${\Cal F}_{b,N} $ where 
${\Cal F}_{b,N} (x) = f_{b,N}(xb^{N+1})$. For $k \leq N/2$, set $X_k = x_k 
b^{-(N+1)}$, so that ${\Cal F}_{b,N} (-X_k) = f_{b,N}(-x_k)$. Hence $\sign 
{{\Cal F}_{b,N} (-X_k) } = (-1)^k$. Next, we note that $0 = X_0 < X_1 < X_2 < 
\dots $; moreover, $X_k = x_1 b^{2k-2 - N -1}$. Hence if $2k < N$, i.e., $2k 
\leq N-1$, then $X_k < 1$ (as $x_1 < b^4$). Now we consider the two cases, 
$N$ even and $N$ odd.

If $N = 2r$ is even, then we have $0 = X_0 < \dots < X_{r-1} < 1$; set $X_r 
= 1$. Then ${\Cal F}_{b,N}(1) $ is up to positive scalar multiple, 
$\sign{(-x)^r) }(1 - 2/b + 2/b^4 -   \dots)$. Now $1 - 2/b + 2/b^4 - 2/b^9 > 
0$ if $b > 1.44$, and since we have assumed $b > \sqrt{3} \equiv 1.73\dots$, 
it follows easily that $\sign{\Cal F_{b,N}(1)} = (-1)^r$. Now for $r + 1 
\leq k \leq N-1 = 2r-1$, set $X_k = (X_{2r-k})^{-1}$. Then $X_j < X_{j+1}$ 
for all $0 \leq j \leq N-2$, and for $k > r$, we have $\sign{\Cal 
F_{b,N}(-X_k)} = (-1)^N \sign {\Cal F_{b,N}(-X_{2r-k})} = (-1)^{2r-k} = 
(-1)^k$. Finally, there exists sufficiently large $X' > X_{2r-1}$ \st 
$\sign{\Cal F_{b,N}(X)} = (-1)^N = 1$; set $X_N = X'$. We thus have $N+1$ 
sign changes in the values of $\sign{\Cal F_{b,N}}$ on the negative real 
numbers, hence ${\Cal F}_{b,N}$ has at least $N$ distinct real roots, and 
thus these must exhaust them. Since $f_{b,N}$ is a reparameterization of 
${\Cal F}_{b,N}(1)$, the same applies to $f_{b,N}$.

If $N = 2r-1$ is odd, then we have $0 = X_0 < X_1 < \dots X_{r-1}  < 1$. For 
$r \leq k \leq 2r-2$, set $X_k = X_{N-k}^{-1}$, and define $X_N$ to be a 
sufficiently large number that $X_N > X_{2r-2}$ and $\sign{\Cal 
F_{b,N}(-X_N)} = (-1)^N$. Then we have $X_j < X_{j+1}$ for $0 \leq j \leq 
N-1$, and moreover, for $N-1 \geq k \geq r$, we have
$\sign{\Cal F_{b,N}(-X_k)} = (-1)^N \sign{\Cal F_{b,N}(-X_{N-k})} = (-1)^{N 
+ N-k} = (-1)^k$. Thus again ${\Cal F}_{b,N}$ has $N+1$ sign changes on the 
negative reals, so ${\Cal F}_{b,N}$, and therefore $f_{b,N}$ has $N $ 
distinct negative roots.

Now define (for fixed $l$) a function of two variables $G(b,Y)$, via $G(b,Y) = f(Yb^4)$, so that
$$\eqalign{
G(b,Y) &= 1 + b^2 Y + b^2 Y^2 + Y^3 + Y^4b^{-4} + \dots\cr
& = 1 + b^2 Y + b^2 Y^2 + Y^3  + \sum_{j=4}^l \frac {Y^j}{b^{j(j+1) - 4j}}\cr
\frac{\partial G}{\partial b}(b,Y) & = 2(Y^2 + Y)b - \sum_{j=2}^{l/2 -1}
\frac{Y^{2j}}{b^{4j^2 - 6j+1}}\(4j^2 -6j +\frac{((2j+1)(2j+2) -8j-4)Y}{b^{4j}}\)\cr
& \qquad - \cases  Y^l C_l & \text{if $l $ is even; $C_l > 0$}\\
\frac{Y^{l-1}}{l(l-1)-4l+5} \((l-1)l -4j+4 - \frac{(l(l+1) -4l)Y}{b^{2l}}\)& \text{if $l $ is odd}\\
\endcases\cr
}$$
It follows easily that if $l \geq 2$, $0 < Y < 1$, and $\sqrt 3 < b < 2$,
then $\frac{\partial G}{\partial b}(b,-Y) < 0$. 

Suppose that for some $\sqrt 3 \leq b_0 < 2$ and $0 < x_0 < b_0^{4}$, we
have $f_{b_0,l} (-x_0) < 0$. Then $G(b_0,-xb_0^{-4}) = 0$ and $0 <
xb_0^{-4} < 1$. By the previous paragraph, it follows that for all $\sqrt
2 > b > b_0$, we have $G(b,x_0 b_0^{-4}) < 0$ and therefore $f_{b,l} (x_0
b^4/b_0^4) \leq 0$. We have thus shown that if $f_{b_0,l} $ hits zero or
less on the interval $(-b_0^4,0)$ and $\sqrt 3 < b_0 < 2$, then if $2 > b
> b_0$, the function $f_{b,l} $ hits zero or less on the interval
$(-b^4,0)$. 

In particular, if $l$ is even, $N \geq 2l$, $\sqrt 3 <b_0 < 2$, and
$f_{b_0,l} $ hits zero or less on the interval $(-b_0^4,0)$, then for all
$b$ with $2 \geq b \geq b_0$, the function $f_{b,N}$ has only real and
simple zeros.

(b) Suppose now that $l$ is odd, and for all $x$ with $0 < x < b^4$ where $b
>\sqrt 3$, $f_{b,l}(-x) \geq 0$. We show that this implies $f_{b,N}$ has
nonreal roots (when $N \geq 2l$), and the same is true when $b$ is
decreased. 

We show that $f'_{b,N}$ (the derivative is \wrt $x$) has a zero in
$(-b^4, 0)$, and that $f_{b,N}$ is strictly positive on $[-b^4,0]$. If
$f_{b,N}$ had only real zeros, this yields a contradiction, since in that
case, the zeros of $f$ are intertwined by those of $f'$. 

We write
$$\eqalign{
f_{b,N} (-x) & = f_{b,l}(-x) +
\sum_{j=1}^{N-l}\frac{(-x)^{l+j}}{b^{(l+j)(l+j+1)}}\cr 
& =f_{b,l}(-x) +
\frac{x^l}{b^{l^2 +l}}\sum_{j=1}^{N-l}\frac{(-1)^{j+1}
x^{j}}{b^{(j)(2l+j+1)}}\cr 
}$$ 
The series is alternating and its
first term ($j=1$) is positive. We check that the series is monotone
decreasing in absolute value; this boils down to
$$\eqalign{
\frac{x^{j}}{b^{j(j+2l+1)}} &
>\frac{x^{j+1}}{b^{(j+1)(j+2l+2)}};\qquad \text{that is,}\cr 
x &<
b^{2j + 2 + l},\cr 
}$$
for which $x < b^4$ is more than sufficient. 
Hence the alternating series is at least as large as the sum of its first two terms, which is
positive. Since $f_{b,l}(-x) \geq 0$ by hypothesis, we have that $f_{b,N}(-x) > 0$. 

Next, we see that $f'_{b,N} (0) = 1/b^2 > 0$ and 
$$\eqalign{
f'_{b,N}(-b^4)& = \sum_{j=1}^N (-1)^{j-1} j b^{-(j^2 - 3j+4)} \cr
& = -b^{-2} + 3b^{-4} - 4b^{-8} + 5  b^{-14} \cr
& < 0 \cr
}$$
as $b^2 > 3$. Hence $f'$ has a zero in the interval $(-b^4,0)$. Since
$f_{b,N}$ is strictly positive on $[-b^4,0]$, it easily follows that
$f_{b,N}$ has nonreal roots.

Next, suppose that $f_{b_0,N}$ is strictly positive on $[-b_0^4,0]$
and $b < b_0$. We show that $f_{b,N}$ is strictly positive on
$[-b^4,0]$. With the same $G$ as defined previously, since
$\frac{\partial G}{\partial b}(b,-Y) < 0$ on the relevant interval,
the result is immediate. 
\qed

The optimal procedure is to look for those values of $b \eqv b_0$ in
the interval $(\sqrt 3,2)$ for which the polynomial (in $x$)
$f_{b,l}$ has a multiple zero at a point in $(0,b^4)$. The even
values of $l$ give lower bounds, the odd values give upper bounds. We
find convergence is extremely fast---by $l=11$, we are within
$10^{-24}$ of $B_0$, the critical value: if $b > B_0$, then $f_{b,N}$
has all of its seros real and simple for all sufficiently large $N$
(including $N = \infty$), and if $b < B_0$, then $f_{b,N}$
has nonreal zeros for all sufficiently large $N$. If $N= \infty$, it
follows that $f_{B_0,N}$ has only real zeros (possibly with
multiplicities), but if $N <\infty$, then it necessarily has nonreal
zeros. It turns out that $B_0$ is just less than $\sqrt{1 + \sqrt
5}$. Remember that the \quotes{$\beta$} value (the constant ratio
$c_j^2/c_{j+1}c_{j-1}$) is the {\it square\/} of $b$. 

To check for multiple roots, we use the discriminant and {\it Maple.}

First for odd $l$; for $l=3$, $b_0 = \sqrt 3$, i.e., $f_{b_0,3}$ has
a multiple (in fact, a triple) zero in the interval $(0,9)$. This
would yield that if $b < \sqrt 3$, then $f_{b,N}$ has nonreal roots
for $N \geq 7$ if we extended the proposition to cover the endpoint
(which is a nuisance). It is easy to show this anyway (keep in mind
however that
$(1+x)^3$ is a reparameterization of $f_{\sqrt 3,3}$, and it has real
zeros, albeit multiple). This computation can be done by hand.

The remaining values were obtained by {\it Maple,} using $50$-digit
accuracy (truncated to $25$ digits here):
$$\eqalign{
l = 5 &\qquad b_0 = 
               1.7982270324863302995970201 \cr
l = 7 &\qquad b_0 =  1.7982315382687507032044628
\cr
l = 9 &\qquad b_0 =  1.7982315382745004887263767
\cr
l = 11 &\qquad b_0 = 1.7982315382745004887933797
\cr
}$$

Now for even $l$; with $l =2$, we obtain $b_0 =2$, which yields
nothing we didn't already know, namely that if all the ratios are
four or more, all the roots are real and simple. However, $l = 4$
yields $b_0 = \sqrt {1 + \sqrt 5}$ (computable by hand), which means
that if $b^2 \geq 1+ \sqrt 5 \eqv (1.79891\dots)^2$ ($3.236\dots$),
then $f_{b,N}$ has only real and simple zeros for $N \geq 8$. Notice
how close this number is to the lower bound obtained from $l=5$ (they
differ at the fourth place). The rest of the values were computed by
{\it Maple,} with the same constraints as above. Here we also have to verify
that $x_0$, one the values of $x$ where the multiple root occurs, can be chosen
in the interval $(0,b_0^4)$. Since $b_0^4 > 10$ (from the lower bounds obtained
in the odd cases), we only need a rough approximation to the values of $x_0$. 

When $l = 4$, $f_{(1+\sqrt 5)^{1/2},4} (x)$ is a quartic and a square, whose roots are
$$
\brcs{-14-6 5^{1/2}+2(50+22 5^{1/2})^{1/2}, -14-6 5^{1/2}-2(50+22 5^{1/2})^{1/2} },
$$
 each with
multiplicity two. The relevant double zero is the first one, so $x_0 = -14-65^{1/2}-2(50+22
5^{1/2})^{1/2}$, approximately $7.49722\dots$ (remember that $x_0$ is the negative of the zero of
the polynomial). This is well within the upper bound of $(1+\sqrt 5)^2$. For $l=6$, the
discriminant of $f_{b,6}$ is up to multiplication by $b^{-210}$, an even polynomial in $b$ of
degree $70$, but even so, we obtain $x_0  \sim 7.503\dots$, and for $l = 8$ and $10$, the
corresponding values of $x_0$ are the same to within six decimals. 
$$\eqalign{
l = 4 &\qquad b_0 = 
               1.7989074399478672722612275 \cr
l = 6 &\qquad b_0 =  1.7982315474312892803918067
\cr
l = 8 &\qquad b_0 =  1.7982315382745016049847445
\cr
l = 10 &\qquad b_0 = 
               1.7982315382745004887933809
\cr
}$$

Now $B_0$ is squeezed between the supremum of the  values of $b_0$ for
odd
$l$ and the infimum of the values of $b_0$ for even $l$. Taking the
numbers from $l = 10$ and $l=11$, we have 
$$
|B_0 -1.7982315382745004887933803| < 6 \times 10^{-25}
$$
{\it Maple\/} (and consequently, I) gave up at $l = 12$, but
presumably it would have yielded accuracy of order $10^{-31}$. In any case,
$B_0^2 \sim 3.23364$ is likely good enough.

As an aside, the self-replicating property discussed above for the polynomials
is more clearly seen in the functional equation satisfied by the entire function $f_{b,\infty}$,
specifically,
$$
f(x) = 1 + \frac x{b^2} f\(\frac x{b^2}\),
$$
which   can be iterated to more confusing forms.

An interesting phenomenon can be observed. From section 1, we can find  polynomials of all
sufficiently large degrees  or an entire function for which all the ratios $c_j^2/c_{j+1}c_{j-1}$
exceed $3.99$, yet which has nonreal zeros. On the other hand, if all these ratios are equal to
$3.24$ (and the degree is large enough), all the zeros are real. This seems
counter-intuitive---the larger the ratios, the better behaved we expect the polynomial to be
\wrt its zeros. However, if we examine the arguments in section 1, we see that when the ratios are large, the coefficients
tail off very quickly, and so only contiguous islands of coefficients with the same signs (see
the argument there) can be used.

\vskip 10pt
\parindent = 1 em

\noindent{\bf References }

\vskip 4pt
\itemitem{[CC]} T Craven \& G Csordas, {\it Complex zero decreasing sequences,} Methods Appl. Anal. 2 (1995),
420--441.

\itemitem{[K]} David C  Kurtz, {\it A sufficient condition for all the roots of a
polynomial to be real,} American Mathematical Monthly, Vol\,99, \#\,3
(1992), 259--263. 

\itemitem{[H]} David Handelman, {\it Isomorphism and non-isomorphsms of AT actions,} J d'Analyse mathŽmatique (to
appear).

\itemitem{[P]} M PŽtrovitch, {\it Une classe remarquable de sŽries entires,} Atti del IV Congresso Internationale
dei Matematici, Rome, Ser 1, 2 (1908), 36--43.

\vskip 15pt

\noindent Mathematics Dept, University of Ottawa, Ottawa ON  K1N 6N5, Canada; e-mail: dehsg\@uottawa.ca

\end